\documentstyle[12pt]{article}
\newcounter{myownsection}
\def\myownsection{\refstepcounter{myownsection} \setcounter{equation}{0}}

\begin{document}
$\;$\\[20pt]
\begin{center}
{\bf NONCOMMUTATIVE EXTENSIONS OF THE \\
FOURIER TRANSFORM AND ITS LOGARITHM} \\[60pt]
{\sc Romuald Lenczewski}\\[40pt]
Institute of Mathematics\\ 
Wroc{\l}aw University of Technology\\
Wybrzeze Wyspianskiego 27\\
50-370 Wroc{\l}aw, Poland\\
e-mail lenczew@im.pwr.wroc.pl\\[40pt]
\end{center}
\begin{abstract}
We introduce noncommutative extensions of the
Fourier transform of probability measures and its logarithm
in the algebra ${\cal A}(S)$ of complex-valued functions on the
free semigroup on two generators $S=FS(\{z,w\})$.
First, to given probability measures $\mu$, $\nu$
whose all moments are finite,
we associate states $\widehat{\mu}$, $\widehat{\nu}$ on 
the unital free *-bialgebra $({\cal B},\epsilon ,\Delta)$ on two
self-adjoint generators $X,X'$ and a projection $P$.
Then we introduce and study cumulants 
which are additive under the convolution 
$\widehat{\mu}\star \widehat{\nu}=\widehat{\mu}\otimes
\widehat{\nu} \circ \Delta$ when restricted
to the ``noncommutative plane'' ${\cal B}_{0}={\bf C}\langle X, X'\rangle$.
We find a combinatorial formula for the M\"{o}bius 
function in the inversion formula and define the 
moment and cumulant generating functions, 
$M_{\widehat{\mu}}\{z,w\}$ and $L_{\widehat{\mu}}\{z,w\}$, 
respectively, as elements of ${\cal A}(S)$.
When restricted to the subsemigroups $FS(\{z\})$ and $FS(\{w\})$, the
function $L_{\widehat{\mu}}\{z,w\}$ coincides with 
the logarithm of the Fourier transform and
with the $K$-transform of $\mu$, respectively. In turn, 
$M_{\widehat{\mu}}\{z,w\}$
is a ``semigroup interpolation'' between the Fourier transform and the
Cauchy transform of $\mu$.
By choosing a suitable weight function $W$ on the semigroup $S$, the moment 
and cumulant generating functions become elements of the Banach algebra
$l^{1}(S,W)$.\\[5pt]
\end{abstract}
\newpage
\myownsection
\begin{center}
{\sc 1. Introduction}
\end{center}
The main examples of noncommutative independence, like tensor, 
free [V1], boolean [Sp-W] and monotone [Mu1] 
lead to different convolutions of measures on 
the real line. This entails existence of 
different cumulants which behave ``nicely'' w.r.t. these convolutions
(for instance, are additive), 
different moment-cumulant formulas and cumulant
generating functions. The latter, like for instance
the $R$-transform and the $S$-transform if free probability [V1,V2]
the $K$-transform for the boolean convolution [S-W] or
the $H$-transform for the monotone convolution [Mu2], 
give noncommutative one-dimensional analogs
of the logarithm of the Fourier transform of probability measures.

However, the connection between them is not so clear -- 
in fact, they arise from quite different theories, 
corresponding to different notions of noncommutative independence 
(for the latter and connections between them, see [G-S], [F], [L1], 
[L3], [S]). Let us mention here that
certain one-parameter interpolations between the combinatorics of
classical and free convolutions have been presented in
[A] and [N]. In turn, two different frameworks including
boolean and free convolutions have also been proposed --
conditional freeness [Bo-Le-Sp] and hierarchy of freeness [L1]
(see also [F-L]). 
On the level of convolutions of measures, 
our motivation, outlined in [L3], can be phrased as follows:
there should exist a noncommutative probability space
with a convolution of states, whose restrictions to certain
one-dimensional commutative subspaces (``real lines'') give 
the known examples of convolutions of measures in 
noncommutative probability.

In this paper the noncommutative space of interest will be the 
unital free *-algebra
$$
{\cal B}_{0}={\bf C}\langle X, X' \rangle 
$$
generated by two self-adjoint generators $X, X'$, on which we will develop
a probability theory. Intuitively, it is helpful to view ${\cal B}_{0}$
as a {\it noncommutative plane}. Using this geometric language,
${\bf C}[X]$ will correspond to the {\it classical real line},
whereas ${\bf C}[X']$ --  to the {\it boolean real line}.
However, in order to make this work, a suitable 
convolution on ${\cal B}_{0}$ has to be introduced.

It is worth noting that most convolutions which appear in noncommutative
probability are, in contrast to the classical convolution, 
highly nonlinear w.r.t. the addition of measures. 
Therefore, in the usual formulation, 
one cannot expect to use *-Hopf algebras or *-bialgebras to define them. 
However, as we showed in [L1] and [L3], 
it is possible to extend a given algebra by a projection $P$ 
and then introduce a *-bialgebra structure on the extended algebra
in many interesting cases. Thus, let
$$
{\cal B}={\bf C}\langle X,X',P\rangle
$$
be the extended free *-bialgebra endowed with the coproduct
$$
\Delta (X)=X \otimes 1 + 1 \otimes X
$$
$$
\Delta (X')=X' \otimes P + P \otimes X', \;\;
\Delta (P)=P\otimes P
$$
and the convolution of states
\begin{equation}
\label{1.1}
\widehat{\mu}\star \widehat{\nu}:= \widehat{\mu}\otimes \widehat{\nu}
\circ \Delta
\end{equation}
called {\it filtered convolution} [L1].

To given measures $\mu , \nu$ on the real line whose all moments
are finite, one can associate states $\widehat{\mu},\widehat{\nu}$ 
on ${\cal B}$ with $P$ playing the role of a ``separator'' of words 
(see Definition 3.1). Then, the restrictions of the convolution (1.1)
to the subalgebras ${\bf C}[X]$ and ${\bf C}[X']$ of ${\cal B}$ 
give classical and boolean convolutions of $\mu$, $\nu$, respectively.
The non-linearity is then ``hidden'' in the definition of 
$\widehat{\mu}, \widehat{\nu}$.
Thus the projection $P$ serves only as a tool to define
appropriate convolutions. One can say that in the usual formulation
of the boolean convolution, where $P$ does not appear, 
one ``sees'' the boolean real line after performing calculations 
involving $P$.

Since ${\cal B}_{0}$ is freely generated
by two generators, it is quite natural to expect 
that our construction should lead to
a more noncommutative cumulant generating function. In fact, 
we show that it is an element of the noncommutative semigroup 
algebra ${\cal A}(S)$ of functions on the
free semigroup $S=FS(2)=FS(\{z,w\})$ on two letters $z,w$,
with the convolution multiplication. We allow in $S$ the empty word
denoted $1$ and denote $S^{+}=S\setminus \{1\}$.
Let us label the variables $X$ and $X'$ from ${\cal B}$ 
with letters $z,w$, respectively, namely $X=X(z)$ and $X'=X(w)$. 
Thus, if $\widehat{\phi}$ is any state on ${\cal B}$, 
the mixed moments of $X(z)$ and $X(w)$  
in the state $\widehat{\phi}$ can be labelled by words of $S$, namely
$$
M_{\widehat{\phi}}(s)=
\widehat{\phi}(X(s_{1})X(s_{2})\ldots X(s_{n})),\;\;{\rm for}\;\;
s=s_{1}s_{2}\ldots s_{n},
$$
where $s_{i}\in \{z,w\}$, $i=1, \ldots , n$, $n\geq 1$, and
we set $M_{\widehat{\phi}}(1)=1$.

The defining recurrence formula for the  
cumulants $L_{\widehat{\phi}}(s)$, where $s\in S^{+}$, 
corresponding to the convolution
(1.1), which we call {\it admissible cumulants}, reads
$$
M_{\widehat{\phi}}(s)=\sum_{p=1}^{l(s)}
\sum_{u=(u_{1}, \ldots , u_{p})\in {\cal AP}(s)}
L_{\widehat{\phi}}(u_{1})\ldots L_{\widehat{\phi}}(u_{p})
$$
where ${\cal AP}(s)$ denotes the set of 
{\it admissible partitions} of the word
$s$, i.e. those which do not have inner $w$'s (see Definition 2.3),
and $l(s)$ denotes the lenght of $s$.

Let $\widehat{\mu}$, $\widehat{\nu}$ be the states on ${\cal B}$ 
associated with $\mu , \nu$ and given by the moments
$$
M_{\widehat{\mu}}(s)=\mu_{n},\;\; M_{\widehat{\nu}}(s)=\nu_{n},
\;\; s=s_{1}\ldots s_{n}\in S
$$
on the noncommutative plane ${\cal B}_{0}$, where
$\mu_{n}, \nu_{n}$ are the $n$-th moments of $\mu ,\nu$, 
respectively, and then extended to ${\cal B}$ by treating $P$ 
as a ``separator'' of words (see Definition 3.1).
Then
\begin{equation}
\label{1.2}
L_{\widehat{\mu}\star\widehat{\nu}}(s)=
L_{\widehat{\mu}}(s) + L_{\widehat{\nu}}(s)
\end{equation}
for any $s \in S^{+}$. Thus, admissible cumulants are 
additive under the convolution (1.1) of states $\widehat{\mu}, \widehat{\nu}$.

For any state $\widehat{\phi}$ on ${\cal B}$, 
the moment and cumulant generating functions are defined as
\begin{eqnarray*}
M_{\widehat{\phi}}\{z,w\}
&=&
\sum_{s\in S}\frac{M_{\widehat{\phi}}(s)}{n(s)!}s,\\
L_{\widehat{\phi}}\{z,w\}
&=&
\sum_{s\in S^{+}}\frac{L_{\widehat{\phi}}(s)}{n(s)!}s,
\end{eqnarray*}
respectively, where 
$$
n(s)!=n_{1}!n_{2}!\ldots n_{p}!,\;\;\; {\rm for}\;\;\; 
s=z^{n_{1}}w^{k_{1}}z^{n_{2}}w^{k_{2}}\ldots w^{k_{p-1}}z^{n_{p}}
$$
with $n_{1},n_{p}\in {\bf N}_{0}:={\bf N}\cup\{0\}$ and 
$k_{1},n_{2}, k_{2}, \ldots , n_{p-1},k_{p-1}\in {\bf N}$.
These formal 
sums should be interpreted as elements of the algebra ${\cal A}(S)$,
a noncommutative ``two-di\-men\-sio\-nal'' analog of the formal power series
${\bf C}[[z]]$.

In this framework, the moments and cumulants corresponding
to the classical and boolean convolutions are labelled 
now not by integers but rather by
elements of the infinite cyclic subsemigroups 
$$
S(z)=FS(\{z\}),\;\;\; S(w)=FS(\{w\})
$$
with $z$ and $w$ as their generators, respectively.
If we restrict the supports of 
$M_{\widehat{\phi}}\{z,w\}$, $L_{\widehat{\phi}}\{z,w\}$ 
to $S(z)$ and $S(w)$, we obtain moment and cumulant generating
functions for the classical and boolean cases, respectively
(the order $n$ of the moments and cumulants in the usual formulation
corresponds to $z^{n}$ and $w^{n}$). In particular, on these restricted
supports, (1.2) gives additivity of the classical and boolean 
cumulants. Thus, one can say that $M_{\widehat{\mu}}\{z,w\}$ is 
a ``semigroup interpolation'' between the Fourier transform 
$F_{\mu}\{z\}=M_{\widehat{\mu}}\{z,0\}$ and the 
Cauchy transform $G_{\mu}(1/w)= wM_{\widehat{\mu}}\{0,w\}$ 
(the left-hand side is treated as a formal power series in $1/w$).
In turn, $L_{\widehat{\mu}}\{z,w\}$ is a ``semigroup interpolation''
between the logarithm of the Fourier transform and 
the $K$-transform of the measure $\mu$.

More generally, one can take a free *-algebra in infinitely many 
indeterminates and then, in this more general framework, define
a similar convolution which unifies, 
apart from tensor and boolean convolutions,
also $m$-free convolutions [F-L] which approximate weakly the free 
convolution (a similar approach includes
the monotone convolution, see [F]).
In the general case we expect to give a ``universal'' noncommutative 
transform of states on a noncommutative version of ${\bf R}^{\infty}$,
whose special cases would also be the $R$-transform and the $H$-transform. 
Nevertheless, it follows from the construction given in 
[L3] that the connection between the classical case and the 
boolean case seems in our approach to be of main importance since
the general model will be obtained by taking copies of ${\cal B}$,
although this step is also non-trivial and will be treated in a subsequent
paper.

Let us note that $L_{\widehat{\phi}}\{z,w\}$, 
our noncommutative extension of the logarithm of the
Fourier transform, is quite different from the cumulant generating 
functions considered so far in noncommutative probability.
In particular, $L_{\widehat{\phi}}\{z,w\}$ and $L_{\widehat{\psi}}\{z,w\}$
do not, in general, commute for 
$\widehat{\phi}\neq \widehat{\psi}$.
Apart from that, it seems to be interesting in its own right 
from the combinatorial point of view. 

In principle, the paper is self-contained, although it is a continuation 
of the study originated in [L3] (the stochastic calculus was developed
in [L2]).
We organized this work as follows. In Section 2 we give basic definitions
on combinatorics of words and we introduce admissible partitions.
A closer look at the filtered convolution for $\widehat{\cal B}$ 
is presented in Section 3. 
In Section 4 we introduce admissible cumulants
and prove that they are additive under the filtered convolution. 
In Section 5 we present the inversion formula for admissible cumulants and 
prove a combinatorial formula for the associated M\"{o}bius function.
In Section 6 we give basic facts on the semigroup algebra ${\cal A}(S)$ 
and the Banach algebra $l^{1}(S,W)$. 
Finally, in Section 7 we derive a formula which expresses the cumulant 
generating function in terms of the moment generating function.
\myownsection
\begin{center}
{\sc 2. Combinatorics on words}
\end{center}
Let $X=\{z,w\}$ be a two-element set. 
A {\it word} on $X$ is a finite sequence
$s=s_{1}s_{2}\ldots s_{n}$, where $s_{i}\in \{z,w\}$ for each 
$i=1, \ldots , n$.
The empty word will be denoted by $1$. The lenght of $s$ will be denoted by
$l(s)$. 
The unital free semigroup $FS(X)=FS(2)$
generated by $X$ is the collection of all words from $X$ made into a semigroup
by the juxtaposition product
$$
(s_{1}\ldots s_{n})(t_{1}\ldots t_{m})=s_{1}\ldots s_{n}t_{1}\ldots t_{m}.
$$
From now on we will understand that $S=FS(2)$ and denote 
$S^{+}=S\setminus \{1\}$.

In the sequel we will use the following terminology and notations:\\[5pt]
(i) $s_{j}\in s$ if $s=s_{1}\ldots s_{n}$, i.e. $s_{j}$ 
is a letter in the word $s$,\\
(ii) $t$ is a {\it subword} of $s=s_{1}\ldots s_{n}\in S$ if it is
a subsequence of the form
\begin{equation}
\label{2.1}
t=s_{i_{1}}\ldots s_{i_{k}}, \;\; 1\leq i_{1}<\ldots < i_{k}\leq n, \; k>0
\end{equation}
or if $t=1$ (the empty word),\\
(iii) for subwords $t,r$ of $s$, 
we write $t<r$ if $t=s_{i_{1}}\ldots s_{i_{k}}$, 
$r=s_{j_{1}}\ldots s_{j_{l}}$
and $i_{1}<j_{1}$; we call $t,r$ {\it disjoint} if the sets
$\{i_{1}, \ldots , i_{k}\}$, $\{j_{1}, \ldots , j_{l}\}$ are disjoint,\\ 
(iv) a subword $t$ of $s$ is a {\it factor} of $s$ if
there exist words $r,r'\in S$ such that $s=rtr'$,\\
(v) if $r,t$ are disjoint subwords of $s$, then $r\cup t$ denotes
the subword of $s$ obtained from $s$ by deleting all letters
which are not in $r$ or $t$; if, in addition, $r\cup t=s$, we
also write $r=s\setminus t$.\\
(vi) if $r,t$ are subwords of $s$, then $r\cap t$ denotes the subword
of $s$ obtained from $s$ be deleting all letters which are not in 
both $r$ and $t$.\\
\indent{\par}
It should be stressed that by subwords we mean
{\it subsequences} of the form (2.1) with encoded information
not only about the letters but also about 
the indices $i_{1}, \ldots , i_{k}$. 
Thus, two distinct subwords
may give the same word. For instance, in the word $zwzw$, there is
only one subword equal to $zwz$, namely $s_{1}s_{2}s_{3}$, 
but there are three subwords
equal to $zw$, namely $s_{1}s_{2}$, $s_{1}s_{4}$ and $s_{3}s_{4}$.
This terminology is borrowed from [Lo] and one should remember
that it is not followed by all authors.
Note also that the subword $r\cup t$ 
is an element of the shuffle of the words $r$ and $t$,
the latter being denoted by $r\circ t$ in [Lo]. One can say that
$r\cup t$ is the only element of the shuffle $r\circ t$, where
all the letters of $r$ and $t$ are ``at the right place''.
For instance, if $r=s_{2}s_{4}$, $t=s_{3}s_{5}$ are subwords
of $s=s_{1}s_{2}s_{3}s_{4}s_{5}$, then $r\cup t=s_{2}s_{3}s_{4}s_{5}$.\\
\indent{\par}
{\sc Definition 2.1.}
By a {\it partition} of the word $s=s_{1}\ldots s_{n}\in S^{+}$ we will
understand any sequence
$$
u=(u_{1}, \ldots , u_{m})
$$
where $u_{1}, \ldots , u_{m}$ are disjoint subwords
of $s$ such that 
$$
s=u_{1}\cup \ldots \cup u_{m}, \;\;{\rm and}\;\;
u_{1}<\ldots < u_{m}
$$
where $1 \leq m \leq n$.
We then write $b(u)=m$, i.e. $b(u)$ denotes the number of
words in the partition $u$.  
We denote by ${\cal P}(s)$ the set of all partitions of the word $s$.
By a {\it factorization} of $s$ we will understand a partition
of $s$ in which every subword $u_{k}$ is a factor, in which case we shall
write $s=u_{1}u_{2}\ldots u_{p}$.\\
\indent{\par}
{\it Remark 1.}
For fixed $s\in S^{+}$ of lenght $n$,
there is a one-to-one correspondence between ${\cal P}(s)$ and 
all partitions ${\cal P}_{n}$ of the set $\{1, \ldots , n\}$
-- subwords correspond to blocks.
Therefore, terminology which refers to partitions of ${\cal P}_{n}$
has its natural analogs in the case of ${\cal P}(s)$. In particular,
we will say that $u\in {\cal P}(s)$ is finer (coarser) than 
$v\in {\cal P}(s)$ if the corresponding partitions of ${\cal P}_{n}$
have this property.\\
\indent{\par}
{\it Remark 2.}
We will adopt the convention that the one-subword 
partition of $s$ consisting of $s$ will be denoted by $s$
instead of $(s)$.\\
\indent{\par}
{\it Remark 3.}
If $u'\in {\cal P}(s')$, $u''\in {\cal P}(s'')$, where
$s',s''$ are disjoint and $s'\cup s''=s$, 
then we will denote by $u'\cup u''$ the 
partition of $s$ consisting of subwords of $s$ which are in 
$u'$ and subwords of $s$ which are in $u''$.
Finally, if $u'$ is a (not necessarily proper) refinement
of $u$, obtained by dividing the subwords of $u$ into perhaps
smaller subwords, we will write $u'\preceq u$.\\
\indent{\par}
{\sc Definition 2.2.}
Let $s=s_{1}\ldots s_{n}\in S^{+}$ and let 
$u=(u_{1}, \ldots , u_{m}) \in {\cal P}(s)$.
We will say that the letter 
$s_{j}$ from the word $u_{p}$ 
is {\it inner} with respect to the word 
$u_{l}=s_{i_{1}}\ldots s_{i_{k}}$, where $l\neq p$, if
$i_{1}<j<i_{k}$. We will then also say that this letter is inner
with respect to the partition $u$.
By a {\it cumulant subword} of $s$ we will understand every subword $r$ 
of $s$ which does not have any inner $w$'s in $s$.
The set of all cumulant subwords of $s$ will be denoted by $C(s)$.
All $w$'s of a cumulant word $r$ will be called $w$-{\it legs}
of $r$. \\ 
\indent{\par}
{\sc Definition 2.3.}
A partition $u=(u_{1}, \ldots , u_{m})\in {\cal P}(s)$, where
$s\in S^{+}$, will be called {\it admissible}, 
if and only if $u_{1}, \ldots , u_{m}$ are cumulant subwords of $s$.
By ${\cal AP}(s)$ we denote the subset of ${\cal P}(s)$ consisting
of admissible partitions. 
If $u\in {\cal P}(s)\setminus {\cal AP}(s)$,
then we will say that there is a {\it non-admissible
inversion} in the sequence of subwords $(u_{i_{1}}, \ldots , u_{i_{p}})$,
where $1\leq i_{1}<i_{2}<\ldots < i_{p}\leq m$, 
if there exists a $w$ in one of these subwords which is inner
w.r.t. another subword. \\
\indent{\par}
{\it Example 1.}
Let $s=z^{2}wz$ and consider two partitions:
$$
u=(s_{1}s_{3},s_{2},s_{4})=(zw,z,z),\;\;\;
u'=(s_{1},s_{2}s_{4},s_{3})=(z,z^{2},w)
$$
given by the diagrams\\
\unitlength=1mm
\special{em:linewidth 0.4pt}
\linethickness{0.4pt}
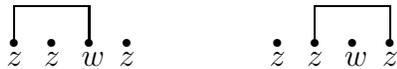
\begin{picture}(120.00,15.00)(-5.00,5.00)
\put(45.00,15.00){\line(1,0){10.00}}
\put(45.00,15.00){\line(0,-1){5.00}}
\put(55.00,15.00){\line(0,-1){5.00}}
\put(45.00,10.00){\circle*{1.00}}
\put(50.00,10.00){\circle*{1.00}}
\put(55.00,10.00){\circle*{1.00}}
\put(60.00,10.00){\circle*{1.00}}
\put(44.00,7.00){$z$}
\put(49.00,7.00){$z$}
\put(54.00,7.00){$w$}
\put(59.00,7.00){$z$}
\put(85.00,15.00){\line(1,0){10.00}}
\put(85.00,15.00){\line(0,-1){5.00}}
\put(95.00,15.00){\line(0,-1){5.00}}
\put(80.00,10.00){\circle*{1.00}}
\put(85.00,10.00){\circle*{1.00}}
\put(90.00,10.00){\circle*{1.00}}
\put(95.00,10.00){\circle*{1.00}}
\put(79.00,7.00){$z$}
\put(84.00,7.00){$z$}
\put(89.00,7.00){$w$}
\put(94.00,7.00){$z$}
\end{picture}
$\;$\\
Then the letter $s_{3}=w$ is inner with respect to
the partition $u'$, but it is not inner w.r.t. $u$.
Therefore, $u$ is admissible, but $u'$ is not admissible.\\
\newpage
\myownsection
\begin{center}
{\sc 3. Moments and convolutions}
\end{center}
Recall from [L1] the definition of the so-called boolean extension of a 
state.\\
\indent{\par}
{\sc Definition 3.1.}
Let $\phi$ be a state on ${\bf C}[Y]$, 
where $Y$ is self-adjoint, i.e. $Y^{*}=Y$. 
The {\it boolean extension} of $\phi$
is the state on ${\bf C}\langle Y, P \rangle $, where $P$ is a projection,
i.e. $P^{*}=P^{2}=P$, given by the linear extension of
\begin{equation}
\label{3.1}
\widetilde{\phi}(P^{\alpha}Y^{n_1}PY^{n_2}P\ldots Y^{n_k}P^{\beta})
=\phi(Y^{n_1})\phi(Y^{n_2})\ldots \phi(Y^{n_k})
\end{equation}
where $\alpha , \beta \in \{0,1\}$ and $n_{1}, \ldots , n_{k}\in {\bf N}$, with
$\widetilde{\phi}(P)=1$. \\
\indent{\par}
The boolean extension of a state is a state since 
$$
\widetilde{\phi}=\phi *_{b}h
$$
($*_{b}$ stands for the boolean product of states),
i.e. $\widetilde{\phi}$ is the boolean product of 
the state $\phi$ on ${\bf C}[Y]$ and the
unital *-homomorphism $h$ on ${\bf C}[P]$ given by $h(P)=h(1)=1$.
That the boolean {\it product} of states is a state, it follows from
the more general case of the so-called conditional product of states
studied in [Bo-Le-Sp]. 
The boolean extensions of states serve as a tool to define 
a new type of convolution of states on the unital free *-algebra
on two generators ${\cal B}_{0}$ (given below in Definition 3.3).\\
\indent{\par}
{\sc Definition 3.2.}
Let ${\cal B}={\bf C}\langle X,X',P\rangle$ be the unital *-algebra of 
polynomials in noncommuting self-adjoint variables 
$X$ and $X'$ and a projection $P$. When endowed with 
the coproduct 
$\Delta:\; {\cal B}\rightarrow {\cal B}\otimes {\cal B}$
and counit $\epsilon: {\cal B}\rightarrow {\bf C}$ given by 
\begin{equation}
\label{3.2}
\Delta (X)= X\otimes 1 + 1 \otimes X, \;\;
\Delta (X')=X'\otimes P + P \otimes X',
\end{equation}
\begin{equation}
\label{3.3}
\Delta (P)=P\otimes P, \;\;\; \epsilon (X)=\epsilon (X')=0, \;\;
\epsilon (P)=1,
\end{equation}
it becomes a unital *-bialgebra, called the {\it filtered bialgebra}.\\
\indent{\par} 
{\sc Definition 3.3.}
Let 
$\eta: \; {\cal B}\rightarrow {\bf C}\langle Y,P \rangle$, 
where $Y$ is self-adjoint and $P$ is a projection,
be the unital *-homomorphism given by
$$
\eta (X)=\eta (X')=Y, \;\; \eta (P)=P, \;\; \eta (1)=1
$$
and let $\phi ,\psi$ be states on ${\bf C}[Y]$. 
Then $\widehat{\phi}=\widetilde{\phi} \circ \eta$, 
$\widehat{\psi}=\widetilde{\psi}\circ \eta$ are states on ${\cal B}$.
Their convolution
\begin{equation}
\label{3.4}
\widehat{\phi}\star \widehat{\psi}=
\widehat{\phi}\otimes \widehat{\psi}\circ \Delta
\end{equation}
will be called the {\it filtered convolution} of $\widehat{\phi}$
and $\widehat{\psi}$.\\
\indent{\par}
We will use the semigroup $S$ to label the mixed moments of
the variables $X$ and $X'$. Namely, we label them by $z$ and $w$, respectively,
to get $X=X(z)$ and $X'=X(w)$. For instance
$$
M_{\widehat{\phi}}(zwz)=\hat{\phi}(XX'X), \;\; 
M_{\widehat{\phi}}(z^{2}w)=\hat{\phi}(X^{2}X'), \;\; 
M_{\widehat{\phi}}(w^{2}zw)=\hat{\phi}((X')^{2}XX').
$$
Note that in our notation $M_{\widehat{\phi}}(s)$ 
(or, simply $M(s)$)
is a moment, not the generating function.
For the latter, we will use the notation $M_{\widehat{\phi}}\{z,w\}$,
or simply $M\{z,w\}$
(see Section 7).

In particular, let now 
$\phi$ and $\psi$ be states on ${\bf C}[Y]$ associated with measures
on the real line $\mu$ and $\nu$, whose all moments are 
finite, i.e.
$$
\phi(Y^{n})=\mu_{n}=\int_{{\bf R}}y^{n}d\mu (y), \;\;\;
\psi(Y^{n})=\nu_{n}= \int_{{\bf R}}y^{n}d\nu(y)
$$
and let $\widetilde{\phi}, \widetilde{\psi}$ be their boolean extensions,
respectively.

Thus, with $\mu$ and $\nu$ one can associate states 
\begin{equation}
\label{3.5}
\widehat{\mu}:=\widehat{\phi}=\widetilde{\phi}\circ \eta,
\;\;\;
\widehat{\nu}:=\widehat{\psi}=\widetilde{\psi}\circ \eta
\end{equation}
on ${\cal B}$ and thus, by restriction, on its subalgebra
${\cal B}_{0}={\bf C}\langle X, X'\rangle$, the unital free 
*-algebra on two generators with moments
$$
M_{\widehat{\mu}}(s)=\widehat{\mu}
(X(s_{1})X(s_{2})\ldots X(s_{n}))=\mu_{n}
$$
$$
M_{\widehat{\nu}}(s)=\widehat{\nu}
(X(s_{1})X(s_{2})\ldots X(s_{n}))=\nu_{n}
$$
for every word $s=s_{1}\ldots s_{n}\in S$, i.e. 
the mixed moments only depend on 
the lenght of $s$ and agree with the moments of the corresponding 
measures on the real line. Informally, one can view them as moments of 
the ``two-dimensional measures'' on the noncommutative plane
${\cal B}_{0}={\bf C}\langle X,X' \rangle $,
canonically associated with measures $\mu$ and $\nu$. 

The filtered convolution $\widehat{\mu}\star \widehat{\nu}$
may also be viewed as a ``two-dimensional measure''
on the noncommutative plane ${\cal B}_{0}$ with moments
\begin{equation}
\label{3.6}
\widehat{\mu} \star \widehat{\nu}(s):
=\widehat{\mu}\star\widehat{\nu} (X(s_{1})X(s_{2})\ldots
X(s_{n})).
\end{equation}
Below we give examples of 
lowest order mixed moments of the filtered convolution
of states $\widehat{\mu}$ and $\widehat{\nu}$.\\
\indent{\par}
{\it Example 1.}
The moments of order $1$ and $2$ do not depend on $s$, namely
\begin{eqnarray*}
\widehat{\mu}\star\widehat{\nu}(s_{1})&=& \mu_{1}+\nu_{1}\\
\widehat{\mu}\star\widehat{\nu}(s_{1}s_{2})&=&\mu_{2}+\nu_{2} +2\mu_{1}\nu_{1}
\end{eqnarray*}
There are two different expressions for moments of order $3$:
\begin{eqnarray*} 
\widehat{\mu}\star\widehat{\nu}(s_{1}zs_{3})&=&\mu_{3}+3\mu_{2}\nu_{1}+
3\mu_{1}\nu_{2} + \nu_{3}\\
\widehat{\mu}\star\widehat{\nu}(s_{1}ws_{3})&=&\mu_{3}+2\mu_{2}\nu_{1}+
\mu_{1}^{2}\nu_{1}+\nu_{3}+2\mu_{1}\nu_{2}+\mu_{1}\nu_{1}^{2}
\end{eqnarray*}
for any $s_{1},s_{3}\in \{z,w\}$,
which coincide with the moments of classical and boolean convolutions
of measures $\mu$, $\nu$,
respectively. In the case of moments of order $4$, we get three
different possibilities:
\begin{eqnarray*}
\widehat{\mu}\star\widehat{\nu}(s_{1}z^{2}s_{4})
&=&
\mu_{4}+4\mu_{3}\nu_{1}+6\mu_{2}\nu_{2}+4\mu_{1}\nu_{3}+\nu_{4}\\
\widehat{\mu}\star\widehat{\nu}(s_{1}w^{2}s_{4})
&=&
\mu_{4}+2\mu_{3}\nu_{1}+2\mu_{2}\nu_{2}+2\mu_{1}\nu_{3}+\nu_{4}\\
&+& 
2\mu_{2}\nu_{1}\mu_{1}+2\mu_{1}\nu_{2}\nu_{1} +\mu_{2}\nu_{1}^{2}
+\nu_{2}\mu_{1}^{2}+2\mu_{1}^{2}\nu_{1}^{2}\\
\widehat{\mu}\star\widehat{\nu}(s_{1}zws_{4})
&=&
\mu_{4}+3\mu_{3}\nu_{1}+2\mu_{2}\nu_{2}+3\mu_{1}\nu_{3}+\nu_{4}\\
&+&
\mu_{1}\nu_{1}\mu_{2}+2\mu_{1}^{2}\nu_{2}+\nu_{1}\mu_{1}\nu_{2}
+2\nu_{1}^{2}\mu_{2}.
\end{eqnarray*}
for any $s_{1},s_{4}\in \{z,w\}$. It is easy to see that
$$
\widehat{\mu}\star\widehat{\nu}(s_{1}wzs_{4})=
\widehat{\mu}\star\widehat{\nu}(s_{1}zws_{4}),
$$
se we have altogether 3 different
cases. The first one corresponds to the classical convolution, the second one 
-- to the boolean convolution, whereas the third one is of a different type.\\
\indent{\par}

Of course, $\widehat{\mu}$, $\widehat{\nu}$ and
$\widehat{\mu}\star \widehat{\nu}$ are defined on all of 
${\cal B}$ and it is not hard to express all their moments in terms of 
the moments on ${\cal B}_{0}$. \\
\indent{\par}
{\sc Proposition 3.4.}
{\it Let  $\widehat{\sigma}\in \{\widehat{\mu}, \widehat{\nu},
\widehat{\mu}\star\widehat{\nu}$\}. Then}
$$
\widehat{\sigma}(P^{\alpha}X(t_{1})PX(t_{2})P\ldots PX(t_{p})P^{\beta})
=\widehat{\sigma}(X(t_{1}))\ldots \widehat{\sigma}(X(t_{p}))
$$
{\it where $\alpha, \beta\in \{0,1\}$ and $t_{i}\in S$, $i=1, \ldots , p$
and the abbreviated notation}
$$
X(t)=X(s_{1})X(s_{2})\ldots X(s_{r})
$$ 
{\it for  $t=s_{1}s_{2}\ldots s_{r}\in S$ is used.}\\[5pt]
{\it Proof.}
This is a straightforward consequence of the definition of $\eta$,
the fact that $P$ acts as a ``separator'' of words in ${\cal B}_{0}$
and that it is group-like, hence the convolution preserves this property.
\hfill $\Box$\\
\indent{\par}
Therefore, we will restrict our attention to the moments of
these states on ${\cal B}_{0}$ since $P$ only serves as a tool to define
a convolution on ${\cal B}_{0}$.\\
\indent{\par}
{\sc Proposition 3.5.} 
{\it Let $\mu$ and $\nu$ be probability measures on the real line, whose
all moments are finite. Then}
$$
\widehat{\mu}\star\widehat{\nu}(s)=
\left\{
\begin{array}{lll}
(\mu\star \nu)_{n} & {\rm if} & s=z^{n}\\
(\mu\uplus \nu)_{n} & {\rm if} & s=w^{n}
\end{array}
\right.
$$
{\it where $n\geq 0$, i.e. the moments of the filtered convolution restricted
to the cyclic subsemigroups $S(z)$ and $S(w)$, 
agree with the moments of the classical and boolean convolutions of
$\mu , \nu$, denoted $\mu\star \nu$ and $\mu\uplus \nu$, respectively.}\\[5pt]
{\it Proof.}
This fact is elementary and follows directly from Definitions 3.1-3.2
(this proposition can also serve as a definition of the boolean convolution).
\hfill $\Box$\\
\myownsection
\begin{center}
{\sc 4. Admissible cumulants}
\end{center}
In this section we define the admissible cumulants and prove that
they are additive under the filtered convolution on ${\cal B}_{0}$.
For notational simplicity, we will denote
the moments and cumulants associated with the state $\widehat{\phi}$
by $M(s)$, $s\in S$, and $L(s)$, $s\in S^{+}$, respectively.\\
\indent{\par}
{\sc Definition 4.1.}
By {\it admissible cumulants} associated with the moments $(M(s))_{s\in S}$
we understand the numbers $(L(s))_{s\in S^{+}}$ defined recursively 
by the formulas
\begin{equation}
\label{4.1}
M(s)=\sum_{p=1}^{l(s)}\sum_{u=(u_{1},\ldots,u_{p})\in {\cal AP}(s)}
L(u_{1})\ldots L(u_{p}),
\end{equation}
where $s \in S^{+}$.\\
\indent{\par}
It is not hard to see that (4.1) is in fact a recurrence formula, which is 
a noncommutative analog of similar recurrence formulas in classical
and noncommutative probability. Namely, we can write
$$
M(s)=L(s)+\sum_{p=2}^{l(s)}\sum_{u=(u_{1}, \ldots , u_{p})\in 
{\cal AP}(s)}L(u_{1})\ldots  L(u_{p})
$$
for any $s \in S^{+}$, and thus $L(s)$ can be expressed in terms
of $M(s)$ and cumulants $L(t)$ associated with words of lenght $l(t)<l(s)$.\\
\indent{\par}
{\it Example 1.} Let $s=s_{1}s_{2}s_{3}s_{4}=zwz^{2}$. Using Definition
4.1 we get
\begin{eqnarray*}
M(zwz^{2})&=&L(zwz^{2})+L(z)L(wz^{2})+L(zw)L(z^{2})+\\
& +& L(zwz)L(z)+ L(z)L(w)L(z^{2})+L(zw)L(z)L(z)+\\
&+& 2L(z)L(wz)L(z)+L(z)L(w)L(z)L(z).
\end{eqnarray*}
Note that there is no contribution to $M(zwz^{2})$ from the partitions
associated with the sequences $(s_{1}s_{4}, s_{2},s_{3})$,
$(s_{1}s_{3}, s_{2}, s_{4})$ and $(s_{1}s_{3}, s_{2}s_{4})$ since in all 
of them $w$ is inner with respect to some block as the figure below
demonstrates:\\ 
\unitlength=1mm
\special{em:linewidth 0.4pt}
\linethickness{0.4pt}
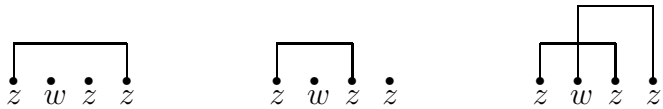
\begin{picture}(120.00,25.00)(-15.00,5.00)
\put(10.00,15.00){\line(1,0){15.00}}
\put(10.00,15.00){\line(0,-1){5.00}}
\put(25.00,15.00){\line(0,-1){5.00}}
\put(10.00,10.00){\circle*{1.00}}
\put(15.00,10.00){\circle*{1.00}}
\put(20.00,10.00){\circle*{1.00}}
\put(25.00,10.00){\circle*{1.00}}
\put(9.00,7.00){$z$}
\put(14.00,7.00){$w$}
\put(19.00,7.00){$z$}
\put(24.00,7.00){$z$}
\put(45.00,15.00){\line(1,0){10.00}}
\put(45.00,15.00){\line(0,-1){5.00}}
\put(55.00,15.00){\line(0,-1){5.00}}
\put(45.00,10.00){\circle*{1.00}}
\put(50.00,10.00){\circle*{1.00}}
\put(55.00,10.00){\circle*{1.00}}
\put(60.00,10.00){\circle*{1.00}}
\put(44.00,7.00){$z$}
\put(49.00,7.00){$w$}
\put(54.00,7.00){$z$}
\put(59.00,7.00){$z$}
\put(80.00,15.00){\line(1,0){10.00}}
\put(85.00,20.00){\line(1,0){10.00}}
\put(80.00,15.00){\line(0,-1){5.00}}
\put(90.00,15.00){\line(0,-1){5.00}}
\put(85.00,20.00){\line(0,-1){10.00}}
\put(95.00,20.00){\line(0,-1){10.00}}
\put(80.00,10.00){\circle*{1.00}}
\put(85.00,10.00){\circle*{1.00}}
\put(90.00,10.00){\circle*{1.00}}
\put(95.00,10.00){\circle*{1.00}}
\put(79.00,7.00){$z$}
\put(84.00,7.00){$w$}
\put(89.00,7.00){$z$}
\put(94.00,7.00){$z$}
\end{picture}
$\;$\\
\indent{\par}
{\it Example 2.} For comparison, take now $s=s_{1}s_{2}s_{3}s_{4}=zw^{2}z$. 
Then
\begin{eqnarray*}
M(zw^{2}z)&=&L(zw^{2}z)+L(z)L(w^{2}z)+L(zw)L(wz) + \\
&+&L(zw^{2})L(z)+L(z)L(w)L(wz)+L(zw)L(w)L(z)+\\
&+& L(z)L(w^{2})L(z)+L(z)L(w)L(w)L(z),
\end{eqnarray*}
and the following partitions give zero contribution:\\
\unitlength=1mm
\special{em:linewidth 0.4pt}
\linethickness{0.4pt}
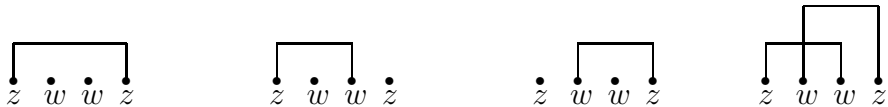
\begin{picture}(120.00,25.00)(-5.00,5.00)
\put(10.00,15.00){\line(1,0){15.00}}
\put(10.00,15.00){\line(0,-1){5.00}}
\put(25.00,15.00){\line(0,-1){5.00}}
\put(10.00,10.00){\circle*{1.00}}
\put(15.00,10.00){\circle*{1.00}}
\put(20.00,10.00){\circle*{1.00}}
\put(25.00,10.00){\circle*{1.00}}
\put(9.00,7.00){$z$}
\put(14.00,7.00){$w$}
\put(19.00,7.00){$w$}
\put(24.00,7.00){$z$}
\put(45.00,15.00){\line(1,0){10.00}}
\put(45.00,15.00){\line(0,-1){5.00}}
\put(55.00,15.00){\line(0,-1){5.00}}
\put(45.00,10.00){\circle*{1.00}}
\put(50.00,10.00){\circle*{1.00}}
\put(55.00,10.00){\circle*{1.00}}
\put(60.00,10.00){\circle*{1.00}}
\put(44.00,7.00){$z$}
\put(49.00,7.00){$w$}
\put(54.00,7.00){$w$}
\put(59.00,7.00){$z$}
\put(85.00,15.00){\line(1,0){10.00}}
\put(85.00,15.00){\line(0,-1){5.00}}
\put(95.00,15.00){\line(0,-1){5.00}}
\put(80.00,10.00){\circle*{1.00}}
\put(85.00,10.00){\circle*{1.00}}
\put(90.00,10.00){\circle*{1.00}}
\put(95.00,10.00){\circle*{1.00}}
\put(79.00,7.00){$z$}
\put(84.00,7.00){$w$}
\put(89.00,7.00){$w$}
\put(94.00,7.00){$z$}
\put(110.00,15.00){\line(1,0){10.00}}
\put(115.00,20.00){\line(1,0){10.00}}
\put(110.00,15.00){\line(0,-1){5.00}}
\put(120.00,15.00){\line(0,-1){5.00}}
\put(115.00,20.00){\line(0,-1){10.00}}
\put(125.00,20.00){\line(0,-1){10.00}}
\put(110.00,10.00){\circle*{1.00}}
\put(115.00,10.00){\circle*{1.00}}
\put(120.00,10.00){\circle*{1.00}}
\put(125.00,10.00){\circle*{1.00}}
\put(109.00,7.00){$z$}
\put(114.00,7.00){$w$}
\put(119.00,7.00){$w$}
\put(124.00,7.00){$z$}
\end{picture}
$\;$\\
\indent{\par}
Thus, among the partitions which do not contribute to
this moment, apart from those of Example 1, we also 
have the partition associated with 
$(s_{1}, s_{2}s_{4},s_{3})$ (the third one in the above figure, it also has 
an inner $w$).\\
\indent{\par}
The restriction of Definition 4.1 to the words from
$S(z)$ gives the usual expression for 
the classical cumulants whereas the restriction to 
$S(w)$ gives the boolean cumulants.
This is beacause ${\cal AP}(z^{n})$ 
can be put in one-to-one correspondence with
all partitions of the set $\{1, \ldots , n\}$, whereas
${\cal AP}(w^{n})$ can be put in one-to-one correspondence
with the interval partitions of $\{1, \ldots , n\}$.
And it is well-known that these two classes of partitions
give classical and boolean cumulants, respectively. 

In the sequel we will need a notation for the summands of 
$\Delta(X(z))$ and $\Delta(X(w))$ given by (3.2)-(3.3):
\begin{equation}
\label{4.2}
j_{1}(X(z))=X(z)\otimes 1, \;\;\; j_{1}(X(w))=X(w)\otimes P
\end{equation}
\begin{equation}
\label{4.3}
j_{2}(X(z))=1\otimes X(z), \;\;\; j_{2}(X(w))=P\otimes X(w)
\end{equation}
(recall that 
$X=X(z)$ and $X'=X(w)$ and compare with the coproduct of Definition 3.2).

Also, for given $s$ of the form (1.4) and given $\epsilon =(\epsilon _{1},
\ldots , \epsilon _{n})$, where $\epsilon _{k}\in \{1,2\}$,
$k=1, \ldots , n$, let
$$
s(1, \epsilon )= \vec{\prod_{j: \epsilon_j =1}}s_{j}, \;\;\;
s(2, \epsilon )= 
\vec{\prod_{j: \epsilon_j =2}}s_{j},
$$
the arrow indicating that the product is taken with the increasing
order of indices. Clearly, $s(1, \epsilon )\cup s(2, \epsilon )$, where
$s(1, \epsilon)$ and $s(2, \epsilon)$ are treated as partitions,
is a partition of $s$.

Below we will give an explicit formula for the mixed moments
\begin{equation}
\label{4.4}
\widehat{\mu}\otimes \widehat{\nu}(s, \epsilon ):=
\widehat{\mu}\otimes \widehat{\nu}
(j_{\epsilon_{1}}(X(s_{1}))\ldots j_{\epsilon_{n}}(X(s_{n})))
\end{equation}
where $\widehat{\mu}, \widehat{\nu}$ are the states 
on ${\cal B}$ given by (3.5) and $\epsilon_{1}, \ldots , 
\epsilon_{n}\in \{1,2\}$.\\
\indent{\par}
{\sc Proposition 4.2.}
{\it Let $\mu , \nu$ be probability measures on the real line with all 
moments finite and let $\epsilon_{1}, \ldots , \epsilon_{n}\in \{1,2\}$. 
Then the mixed moments (4.4) are given by the formula}
\begin{equation}
\label{4.5}
\widehat{\mu}\otimes \widehat{\nu}(s, \epsilon )
=
\widehat{\mu}(u') 
\widehat{\nu}(u'' )
\end{equation}
{\it where }
\begin{eqnarray}
\label{4.6}
\widehat{\mu} (u ')&=&\widehat{\mu} (u_{1}')\ldots \widehat{\mu} (u_{p}')\\
\label{4.7}
\widehat{\nu} (u '' )&=&\widehat{\nu} (u_{1}'')\ldots \widehat{\nu} (u_{q}'')
\end{eqnarray}
{\it and 
$u' =(u_{1}', \ldots , u_{p}')$, 
$u'' =(u_{1}'', \ldots , u_{q}'')$ are the unique
coarsest partitions of $s(1, \epsilon )$ and 
$s(2, \epsilon )$, respectively, which define
an admissible partition of $s$ (their dependence on $\epsilon$ is
suppressed).}\\[5pt]
{\it Proof.}
By substituting (4.2)-(4.3) into (4.4) and using Proposition 3.4, we get
$$
\widehat{\mu}\otimes \widehat{\nu}(s, \epsilon )
=\widehat{\mu}\circ \tau ' (s) \; 
\times 
\widehat{\nu}\circ \tau '' (s)
$$
where 
$$
\tau '(s)=\tau '(s_{1})\ldots \tau '(s_{n}), \;\;\;
\tau ''(s)=\tau ''(s_{1})\ldots \tau ''(s_{n}),
$$
and
$$
\tau '(s_{j})=\left\{
\begin{array}{lll} 
X(s_{j}) & {\rm if} & \epsilon_{j}=1\\
P & {\rm if} & (s_{j}, \epsilon_{j})=(w, 2)\\
1 & {\rm if} & (s_{j}, \epsilon_{j})=(z, 2)
\end{array}
\right.
\;\;\;\;
\tau ''(s_{j})=\left\{
\begin{array}{lll} 
X(s_{j}) & {\rm if} & \epsilon_{j}=2\\
P & {\rm if} & (s_{j}, \epsilon_{j})=(w, 1)\\
1 & {\rm if} & (s_{j}, \epsilon_{j})=(z, 1)
\end{array}
\right.
$$
and therefore the $P$'s define interval 
partitions of $s(1, \epsilon )$ and
$s(2, \epsilon )$, respectively, denoted 
$u' =(u_{1}', \ldots , u_{p}')$ and 
$u'' =(u_{1}'', \ldots , u_{q}'')$, 
where $1 \leq p+q \leq n$. 
The words of these partitions are the longest
subwords of $s(1, \epsilon )$ and $s(2, \epsilon )$ for which
the corresponding products of $X(s_{j})$'s are not separated by
a $P$. The pair $(u',u'')$
defines an admissible partition $u=u'\cup u''=(u_{1}, \ldots , u_{p+q})$ 
of $s$. In fact, each of its subwords, say $u_{k}$, belongs to either 
$u'$ or $u''$ -- without loss of generality we can 
suppose that $u_{k}=u_{r}'$ for some $r$.
Then, between the letters of $u_{k}$, say $s_{j}$ and $s_{l}$, 
there can only be letters of the same block $u_{k}$ or letters of the blocks 
of $u''$. The latter have to be $z$'s
since any $w$ would produce a $P$ between $X(s_{j})$ and $X(s_{l})$
at the first tensor site as the mapping $\tau '$ indicates, but then
$s_{j}$ and $s_{l}$ would not belong to the same block, which is a 
contradiction. This completes the proof. \hfill $\Box$\\
\indent{\par}
{\sc Proposition 4.3.}
{\it Under the assumptions of Proposition 4.2, the mixed moments
$\widehat{\mu}\otimes \widehat{\nu}(s, \epsilon )$ 
can be expressed in terms of cumulants as follows:}
\begin{equation}
\label{4.8}
\widehat{\mu}\otimes \widehat{\nu}(s, \epsilon )=
\sum_{p=1}^{n}
\sum_{u=(u_{1}, \ldots , u_{p})\in {\cal AP}_{\epsilon}(s)}
L_{\epsilon (1)}(u_{1})\ldots L_{\epsilon (p)}(u_{p})
\end{equation}
{\it where}
$$
L_{i}=
\left\{
\begin{array}{lll}
L_{\widehat{\mu}} &{\rm if} & i=1\\
L_{\widehat{\nu}} &{\rm if} & i=2
\end{array}
\right.
$$
{\it and ${\cal AP}_{\epsilon}(s)$ denotes the set of all admissible 
partitions of $s$ which are subpartitions of the partition $(s(1, \epsilon),
s(2, \epsilon ))$ and $\epsilon (k)=i$, $i\in \{1,2\}$,
if for all $s_{j}\in u_{k}$ we have $\epsilon_{j}=i$.}\\[5pt]
{\it Proof.}
By applying Definition 4.1 to every moment $\widehat{\mu}(u_{k}')$
and $\widehat{\nu}(u_{l}'')$ on the RHS of (4.5), i.e. expressing these
moments in terms of cumulants, we obtain 
$\widehat{\mu}\otimes \widehat{\nu}(s, \epsilon )$
equal to a sum of products of type
\begin{equation}
\label{4.9}
L_{\epsilon (1)}(v_{1})\ldots L_{\epsilon (p)}(v_{m}),
\end{equation}
where $v=(v_{1}, \ldots , v_{m}) \in {\cal AP}_{\epsilon}(s)$
is a refinement of $u=u'\cup u''\in {\cal AP}_{\epsilon}(s)$.
Moreover, this refinement must be admissible by
the definition of admissible cumulants (it should be remembered that
a refinement of an admissible partition does not have to be admissible).

Since all products of type (4.9) which are obtained in this fashion
are associated with {\it different} admissible refinements of the partition
$u=u'\cup u''$, they give distinct elements of ${\cal AP}_{\epsilon}(s)$.
Therefore, we just need to prove that on the RHS of (4.5)
we obtain products of cumulants
associated with {\it all} $u\in {\cal AP}_{\epsilon }(s)$.
Thus, let $v=(v_{1}, \ldots , v_{m})\in {\cal AP}_{\epsilon} (s)$.
Take the partition $u'\cup u''$ 
of Proposition 4.2.
It is enough to show that letters (understood as pairs
$(s_{k},k)$) from each word $v_{j}$, $1\leq j \leq m$, 
with say $\epsilon (j)=1$, cannot belong to different words 
of $u'$. Suppose that two letters, say $s_{i},s_{k}\in v_{j}$, belong
to different words of $u'$. This means that they must be separated in $s$
by a $w=s_{l}$ with $\epsilon_{l}=2$ and, therefore, that 
$v$ is not admissible, which is a contradiction. We conclude that
$v$ must be an admissible refinement of $u$. However,
all admissible refinements of $u$ give a contribution to the
RHS of 4.5, hence this ends the proof. \hfill $\Box$\\
\indent{\par}
{\it Example 3.}
Let us give two examples of moments 
$\widehat{\mu}\otimes \widehat{\nu}(s, \epsilon )$ expressed
in terms of cumulants according to (4.8). For the sake of generality,
we take $s=s_{1}s_{2}s_{3}s_{4}$, i.e. an arbitrary word of lenght $l(s)=4$.
We shall use the notation
\begin{equation}
\label{4.10}
\delta_{k}=
\left\{
\begin{array}{lll}
1 & {\rm if} & s_{k}=z\\
0 & {\rm if} & s_{k}=w
\end{array}
\right.
\end{equation}
and, for simplicity, we shall write $L(i_{1}\ldots i_{n})$
instead of $L(s_{i_1}\ldots s_{i_n})$.
Take, for instance, the two most interesting examples of
$\epsilon =(1,1,2,1)$ and $\epsilon '(1,2,1,2)$. We have
\begin{eqnarray*}
\widehat{\mu}\otimes \widehat{\nu}(s, \epsilon )&=&
L_{1}(12)L_{2}(3)L_{1}(4)+ 
\delta_{3}L_{1}(124)L_{2}(3)\\
&+&
\delta_{2}\delta_{3} L_{1}(14)L_{1}(2)L_{2}(3)
+\delta_{3}L_{1}(1)L_{1}(24)L_{2}(3)\\
&+&
L_{1}(1)L_{1}(2)L_{2}(3)L_{1}(4)\\
\widehat{\mu}\otimes \widehat{\nu}(s, \epsilon ')&=&
\delta_{2}\delta_{3}L_{1}(13)L_{2}(24)
+ \delta_{2}L_{1}(13)L_{2}(2)L_{2}(4)\\
&+&
\delta_{3}L_{1}(1)L_{2}(24)L_{1}(3)
+L_{1}(1)L_{2}(2)L_{1}(3)L_{2}(4).
\end{eqnarray*}
In the special cases of $s=z^{4}$ (all $\delta$'s are 
equal to 1) and $s=w^{4}$ (all $\delta$'s vanish), we 
get mixed moments of classical and boolean variables, respectively.
\\
\indent{\par}
Let us show now that the admissible cumulants are additive under the filtered
convolution.\\
\indent{\par}
{\sc Theorem 4.4.} ({\sc Additivity of cumulants})
{\it Let $\mu , \nu$ be probability measures on the real line with all 
moments finite. Let $\widehat{\mu}, \widehat{\nu}$ be the associated states
on ${\cal B}$ given by (3.5). Then}
\begin{equation}
\label{4.11}
L_{\widehat{\mu}\star\widehat{\nu}}(s)
=L_{\widehat{\mu}}(s)+ L_{\widehat{\nu}}(s)
\end{equation}
{\it for every $s\in S^{+}$.}\\[5pt]
{\it Proof.}
We will use the induction argument with respect to the lenght of $s$.
It is clear that if $l(s)=1$, then
$$
L_{\widehat{\mu}\star\widehat{\nu}}(s)
=M_{\widehat{\mu}\star\widehat{\nu}}(s)
=M_{\widehat{\mu}}(s) + M_{\widehat{\nu}}(s)
=L_{\widehat{\mu}}(s) + L_{\widehat{\nu}}(s).
$$
Suppose now that (4.11) holds for words $s$ of lenght $l(s)\leq n-1$. We will show 
that then (4.11) holds for words $s$ of lenght $l(s)=n$. 
By Definition 4.1, we have
$$
L_{\widehat{\mu}\star\widehat{\nu}}(s)=
M_{\widehat{\mu}\star\widehat{\nu}}(s)
-
\sum_{p=2}^{n}\sum_{u=(u_{1}, \ldots , u_{p})\in {\cal AP}(s)}
L_{\widehat{\mu}\star\widehat{\nu}}(u_{1})
\ldots 
L_{\widehat{\mu}\star\widehat{\nu}}(u_{p})
$$
for $s=s_{1}\ldots s_{n}$.

We know from Section 3 that
\begin{eqnarray*}
M_{\widehat{\mu}\star\widehat{\nu}}(s)
&=&\widehat{\mu}\otimes \widehat{\nu}
(\Delta X(s_{1})\ldots \Delta X(s_{n}))\\
&=&
\sum_{\epsilon_{1}, \ldots , \epsilon_{n}\in \{1,2\}}
\widehat{\mu}\otimes \widehat{\nu}
(s, (\epsilon_{1}, \ldots , \epsilon_{n} ))
\end{eqnarray*}
Using the inductive assumption, we have
$$
L_{\widehat{\mu}\star\widehat{\nu}}(u_{k})
=
L_{\widehat{\mu}}(u_{k})+L_{\widehat{\nu}}(u_{k}), \;\;\;
\forall \; k=1, \ldots , p
$$
since $l(u_{k})<n$ for all $1\leq k \leq n$ (recall that $p\geq 2$).
Therefore,
\begin{eqnarray*}
L_{\widehat{\mu}\star\widehat{\nu}}(s)&=&
M_{\widehat{\mu}}(s)+M_{\widehat{\nu}}(s)-
\sum_{p=2}^{n}\sum_{u=(u_{1}, \ldots , u_{p})\in {\cal AP}(s)}
L_{\widehat{\mu}}(u_{1})\ldots L_{\widehat{\mu}}(u_{p})\\
&-&
\sum_{p=2}^{n}\sum_{u=(u_{1}, \ldots , u_{p})\in {\cal AP}(s)}
L_{\widehat{\nu}}(u_{1})\ldots L_{\widehat{\nu}}(u_{p}) +D(s)\\
&=& L_{\widehat{\mu}}(s) + L_{\widehat{\nu}}(s) + D(s)
\end{eqnarray*}
where
\begin{eqnarray*}
D(s)&=&
\sum_{\stackrel{\epsilon_{1}, \ldots , \epsilon_{n}\in \{1,2\}}
{\scriptscriptstyle {\rm not}\;{\rm all}\;{\rm equal}}}
\widehat{\mu}\otimes \widehat{\nu}
(s, (\epsilon_{1}, \ldots , \epsilon_{n}) )
\\
&-&
\sum_{p=2}^{n}\sum_{u=(u_{1}, \ldots , u_{p})\in {\cal AP}(s)}
\sum_{\stackrel{\epsilon (1), \ldots , \epsilon (p)\in \{1,2\}}
{\scriptscriptstyle {\rm not}\;{\rm all}\;{\rm equal}}}
L_{\epsilon (1)}(u_{1})\ldots L_{\epsilon (p)}(u_{p})
\end{eqnarray*}
where $L_{1}(u)=L_{\widehat{\mu}}(u)$ and $L_{2}(u)=L_{\widehat{\nu}}(u)$.

Using Proposition 4.3 and interchanging the summations, which in 
this case takes the form
$$
\sum_{\stackrel{\epsilon_{1}, \ldots , \epsilon_{n}\in \{1,2\}}
{\scriptscriptstyle {\rm not}\;{\rm all}\;{\rm equal}}}
\sum_{u=(u_{1}, \ldots , u_{p})\in {\cal AP}_{\epsilon}(s)}=
\sum_{u=(u_{1}, \ldots , u_{p})\in {\cal AP}(s)}
\sum_{\stackrel{\epsilon (1), \ldots , \epsilon (p)\in \{1,2\}}
{\scriptscriptstyle {\rm not}\;{\rm all}\;{\rm equal}}},
$$
for every $p\geq 2$, we deduce that $D(s)=0$, which completes the proof.
\hfill $\Box$\\
\myownsection
\begin{center}
{\sc 5. M\"{o}bius Inversion Formula}
\end{center}
In this section we apply the theory of M\"{o}bius
functions to prove an inversion formula for the admissible 
cumulants. We also derive a combinatorial formula for 
the associated M\"{o}bius function. 
For details on the theory of M\"{o}bius functions, see [R1] and [R2].\\
\indent{\par}
{\sc Proposition 5.1.}
{\it For every $s\in S^{+}$, the set of admissible partitions
${\cal AP}(s)$ is a lattice.}\\[5pt]
{\it Proof.}
We will show that if $u,u'\in {\cal AP}(s)$, then $u\wedge u' ,
u\vee u' \in {\cal AP}(s)$, where $\wedge $ and $\vee $ denote 
meet and join in the lattice of all partitions of $s$, ${\cal P}(s)$.

Let $s=s_{1}\ldots s_{n}$,
$u=(u_{1}, \ldots , u_{r})$, $u'=(u_{1}', \ldots , u_{r'}')$ 
and let 
$$
u \wedge u'=(v_{1}, \ldots , v_{p}), \;\;
u \vee u' =(t_{1}, \ldots , t_{q}).
$$
We have
$v_{j}=u_{k}\cap u_{l}'$ for some $k,l$. Suppose that there exists
$w=s_{m}\in v_{j}$ which is inner w.r.t. $v_{j'}=u_{k'}\cap u_{l'}$
for $j\neq j'$. But then $s_{m}$ would be inner w.r.t. $u$ or $u'$
since we must have $(k,l)\neq (k',l')$, which would imply that
either $u$ or $u'$ is not admissible, which is a contradiction.

Suppose now that there exists a $w=s_{m}\in t_{j}$ which is inner w.r.t. 
$t_{j'}$, where $j\neq j'$. Then $s_{m}\in u_{k}$
for some $k$. Clearly, there do not exist $s_{r}, s_{r'}\in u_{l}$
where $l\neq k$ and $r<m <r'$ since in that case $s_{m}$ would be inner
w.r.t. $u_{l}$ and $u$ would not be admissible. Therefore, there must exist
$s_{r}\in u_{l}$ and $s_{r'}\in u_{l'}$, with $l\neq l'$ (of course,
$l,l'\neq k$). To fix attention, let $1\leq r<m<r'\leq n$. 
Note that all letters of $u_{l'}$
must follow $s_{m}$ and all letters of $u_{l}$ must precede $s_{m}$
in $s$ by the argument above. 
In a similar manner we can show that every subword
of $u'$ must either precede or follow $s_{m}$.
This implies that the partition of $s$ obtained from $u\vee u'$
by splitting $t_{j'}$ into 
$t_{j'}\cap s_{1}\ldots s_{m-1}$ and $t_{j'}\cap s_{m+1}\ldots s_{n}$
is finer than $u\vee u'$ and coarser than both $u$ and $u'$, which 
is a contradiction. This completes the proof.
\hfill $\Box$\\
\indent{\par}
Each lattice ${\cal AP}(s)$, where $s=s_{1}\ldots s_{n}$, 
has the unique minimal element $0_{s}=(s_{1}, \ldots , s_{n})$
and the unique maximal element $1_{s}=(s)=s$. 
We will often skip the index $s$ in the first notation
if it is clear which lattice is considered,
whereas $s$ will be used instead of $1_{s}$.

In order to apply the theory of M\"{o}bius functions to the combinatorics
of moments and cumulants on ${\cal B}_{0}$, we need to take the union
of the lattices of admissible partitions of all nonempty words, namely
$$
P:=\bigcup_{s\in S^{+}}{\cal AP}(s),
$$
on which we introduce partial order by the condition
$$
u\leq v \;\;{\rm iff}\;\; (\exists \;s\in S^{+}: \;u,v\in {\cal AP}(s)
\;\;{\rm and}\;\; u\preceq v)
$$
where $\preceq $ is the usual partial order 
inherited from ${\cal P}_{n}$ for $l(s)=n$.
As usual, we will write $u<v $ if $u \leq v$ and $u\neq v$.
By the segment $[u,v]$, where $u,v\in {\cal AP}(s)$,
we denote the set of all partitions $t$ such that $u\leq t \leq v$.

By the {\it incidence algebra} of the partially order set
$P$, denoted $I(P)$, 
we will understand the set of complex-valued functions
$$
f:\;  P\times P \rightarrow {\bf C}
$$
with values denoted $f(u|v)$, 
such that $f(u|v)=0$ unless $u\leq v$.
If the second argument of functions from
the incidence algebra $I(P)$ 
is a one-word partition of $s$, then we will often skip the second argument
and write 
$$
f(u|s)=f(u)=f(u_{1}, \ldots , u_{p})
$$
for $u=(u_{1}, \ldots , u_{p})$. \\
\indent{\par}
{\it Example 1.}
Note that the order relation in $P$ is stronger than
taking a refinement. For instance,
$$
u=(s_{1}s_{3},s_{2}s_{5},s_{4})\preceq (s_{1}s_{3}s_{4},s_{2}s_{5})=v
$$
for any $s=s_{1}s_{2}s_{3}s_{4}s_{5}$, but 
$$
u\leq v \;\; {\rm iff} \;\; s_{2}=s_{3}=s_{4}=z
$$
and thus $f(u|v)=0$ unless $s_{2}=s_{3}=s_{4}=z$, for any
$f\in I(P)$.\\
\indent{\par}
The sum and multiplication
by scalars in $I(P)$ are defined as usual. The product is given by
$$
h(u|v)=\sum_{u\leq t \leq v}f(u|t)g(t|v)
$$
and the identity element of the algebra is given by the Kronecker delta
$\delta (u|v)$. The {\it zeta function} of $P$ is defined as
$$
\zeta (u|v)= \left\{
\begin{array}{ll}
1 & {\rm if}\;\; u \leq v\\
0 & {\rm otherwise}
\end{array}
\right.
$$
and the function $i(u|v)=\zeta (u|v)-\delta (u|v)$ is called
the {\it incidence function}.

It is well-known that the zeta function is invertible in the incidence 
algebra. The inverse is called the {\it M\"{o}bius function} and is given by
the recursion
\begin{equation}
\label{5.1}
m(u|v)=\left\{
\begin{array}{cl}
1 & {\rm if}\;\; u=v\\
-\sum\limits_{u \leq t \leq v} m(u|t) & {\rm otherwise}.
\end{array}
\right.
\end{equation}
and the {\it M\"{o}bius inversion formula} reads:
\begin{equation}
\label{5.2}
g(u)=\sum_{v\leq u}f(v) \Longrightarrow f(u)=\sum_{v\leq u}m(v|u)g(v)
\end{equation}
for functions $f,g: P\rightarrow {\bf C}$.

In order to apply the theory of M\"{o}bius functions to invert 
formula (4.1), let us define multiplicative functions 
$M(u)$ and $L(u)$ for $u$ ranging over ${\cal AP}(s)$:
\begin{eqnarray}
\label{5.3}
M(u)&=&M(u_{1})\ldots M(u_{p})\\
\label{5.4}
L(u)&=&L(u_{1})\ldots L(u_{p})
\end{eqnarray}
for $u=(u_{1}, \ldots , u_{p})\in {\cal AP}(s)$. For the applications of 
the M\"{o}bius inversion formula to free probability, see [Sp].\\
\indent{\par}
{\sc Proposition 5.2.}
{\it The partition-dependent moments and cumulants
satisfy the relation}
\begin{equation}
\label{5.5}
M(u)=\sum_{v\leq u}L(v)
\end{equation}
{\it where $M(u)$ and $L(v)$ are given by (5.3)-(5.4)
and $u\in {\cal AP}(s)$, $s\in S^{+}$.}\\[5pt]
{\it Proof.}
Clearly, formula (5.5) holds for $u=s$ by Definition 4.1.
We need to justify that it holds if $u<s$. Using (5.3) and
then expressing every $M(u_{k})$ on the RHS of (5.3) in terms of
cumulants according to (4.1), we get
\begin{equation}
\label{5.6}
M(u)=\sum_{v^{1}\leq u_{1}}\ldots \sum_{v^{p}\leq u_{p}}
L(v^{1})\ldots L(v^{p})
\end{equation}
where $v^{1}, \ldots , v^{p}$ are partitions of $u_{1}, \ldots , u_{p}$, 
respectively,
and we only need to show that the RHS of equation (5.6) can be written
in the form given by equation (5.5). 
First, let us show that
if $v^{1}\leq u_{1}, \ldots v^{p}\leq u_{p}$, then
$v:=v^{1}\cup \ldots \cup v^{p}\in {\cal AP}(s)$.
Recall that $v^{k}\leq u_{k}$ means that
$v^{k}$ is an admissible partition of $u_{k}$, $k=1, \ldots , p$.
Suppose there exists a $w\in v_{j}^{k}$, where $v_{j}^{k}$ is the $j$-th
subword of $v^{k}$, which is inner w.r.t.
$v_{j'}^{k'}$, where $v_{j'}^{k'}$ is the $j'$-th subword of $v^{k'}$.
We must have $k=k'$ since otherwise this $w$ would be inner w.r.t.
the subword $u_{k'}$ which would imply that $u$ is not admissible.
Thus, assume that $k=k'$ (of course, in that case we must have $j\neq j'$).
But then $v^{k}$ is not an admissible partition of $u_{k}$, which
is a contradiction. Therefore, $v\leq u$.

Suppose now that there exists $v<u$ which is not obtained
from $u$ by taking admissible partitions of the words $u_{1}, \ldots , u_{k}$, 
respectively.
Of course, $v=v^{1}\cup \ldots \cup v^{p}$, where
$v^{k}\in {\cal P}(u_{k})$ for $k=1, \ldots , p$.
Since $v\in {\cal AP}(s)$, there is no $w$ in one subword, say
$v_{j}^{k}$, inner w.r.t. another subword $v_{j'}^{k'}$, 
where $(j,k)\neq (j',k')$.
This implies that $v^{k}\in {\cal AP}(u_{k})$.\hfill $\Box$\\
\indent{\par}
More generally, one can show that we have the formula
\begin{equation}
\label{5.7}
[v,u]\cong [v^{1},u_{1}]\times \ldots \times [v^{p},u_{p}]
\end{equation}
where $u=(u_{1}, \ldots , u_{p})$, $v^{k}=v\cap u_{k}$
is the partition of $u_{k}$ consisting of those subwords
of $v$ whose union gives $u_{k}$,
$[v,u]$ is the segment in the lattice ${\cal AP}(s)$, 
$[v^{k},u_{k}]$ is the segment in ${\cal AP}(u_{k})$ with
$u_{k}$ being treated as a subword of $s$.\\
\indent{\par}
{\sc Theorem 5.3.} ({\sc Inversion Formula for Cumulants})
{\it Let $(M(s))_{s\in S}$ be 
the mixed moments on the noncommutative plane 
$\widehat{\cal B}_{0}$ in some state 
$\widehat{\phi}$. Then the corresponding admissible 
cumulants $(L(s))_{s\in S^{+}}$ are given by}
\begin{equation}
\label{5.8}
L(s)=\sum_{u\leq s}m(u)M(u)
\end{equation}
{\it where $M(u)$ is given by (5.3) and $m(u)=m(u|s)$.}\\[5pt]
{\it Proof.}
It is a special case of the general M\"{o}bius inversion formula given
by (5.2), which can be used in view of Proposition 5.2.
\hfill $\Box$\\
\indent{\par}
In the examples given below we compute certain $m(u)$, $u\in {\cal AP}(s)$,
using the formula
\begin{equation}
\label{5.9}
m(u|v)=\zeta^{-1}(u|v)=\delta (u|v) - i(u|v) + i^{2}(u|v) \ldots,
\end{equation}
which expresses the M\"{o}bius function in terms of the incidence
function [R1].\\
\indent{\par}
For notational simplicity we identify $u_{j}$ with number
$j$ and thus use a short-hand notation for words  
$j\cup k=u_{j}\cup u_{k}$, $j\cup k\cup l=u_{j}\cup u_{k} \cup u_{l}$, 
etc. For instance
\begin{eqnarray*}
m(1,2,3)
&=&
m((u_{1},u_{2},u_{3})|u_{1}\cup u_{2}\cup u_{3}),\\
m(1\cup 3,2)
&=&
m((u_{1}\cup u_{3},u_{2})|u_{1}\cup u_{2}\cup u_{3}),
\end{eqnarray*}
where we also skip some parentheses for notational sipmlicity.\\
\indent{\par}
{\it Example 1.}
Let $u=(u_{1},u_{2},u_{3})\in {\cal AP}(s)$.
Then
\begin{eqnarray*}
m(1,2,3)
&=& 
- i(1,2,3)
+
i(1,2,3|1\cup 2,3)i(1\cup 2,3)\\
&&+
i(1,2,3|1\cup 3,2)i(1\cup 3,2)
+
i(1,2,3|1,2\cup 3)i(1,2\cup 3)
\end{eqnarray*}
and, in particular, if the subwords of $u$ are one-letter words, i.e.
$u_{k}=s_{k}$, $k=1,2,3$, this gives
$$
m(1,2,3)=-1+1+\delta_{2}+1=1+\delta_{2}=
\left\{
\begin{array}{lll}
1 & {\rm if} & s_{2}=w\\
2 & {\rm if} & s_{2}=z
\end{array}
\right.
$$
for $s=s_{1}s_{2}s_{3}$, where the notation (4.10) is used.\\
\indent{\par}
{\it Example 2.}
In a similar manner, if blocks are one-letter words, we get 
\begin{eqnarray*}
m(1\cup 3,2,4)
&=& 
-i(1\cup 3,2,4)\\
&&+
i(1\cup 3,2,4|1\cup 2\cup 3,4)
i(1\cup 2\cup 3,4)\\
&&+ 
i(1\cup 3,2,4|1\cup 3\cup 4,2)
i(1\cup 3\cup 4,2) \\
&&+
i(1\cup 3,2,4|1\cup 3,2\cup 4)
i(1\cup 3,2\cup 4)\\
&=& -\delta_{2}+\delta_{2}+\delta_{2}+\delta_{2}\delta_{3}\\
&=&\delta_{2}+ \delta_{2}\delta_{3}\\
&=&
\left\{
\begin{array}{lll}
0 &{\rm if} & s=s_{1}ws_{3}s_{4}\\
1 &{\rm if} & s=s_{1}zws_{4}\\
2 &{\rm if} & s=s_{1}z^{2}s_{4}
\end{array}
\right.
\end{eqnarray*}
where $s_{1},s_{3},s_{4}\in \{z,w\}$.\\
\indent{\par}
Looking at these examples, it is not hard to observe
that when computing the M\"{o}bius function from formula
(5.9), we get cancellations (the number of these gets quite large
when $b(u)$ increases). Therefore, it is important to obtain a simpler 
formula for the M\"{o}bius function, in which all cancellations would
be taken into account. 
This phenomenon is a rather typical but also non-trivial
part of the theory (for some classical examples, see [R1] and [R2]). 
In order to do that, we will introduce the notion of
{\it admissible shuffles} of $u\in {\cal AP}(s)$.
It will turn out that 
$$
m(u)=(-1)^{b(u)-1}a(u),
$$ 
where $a(u)$ denotes the number of admissible shuffles of $u$,
a noncommutative analog of $(p-1)!$ -- 
the number of ways we can shuffle the $p-1$ blocks of a partition
consisting of $p$ blocks, keeping the first block fixed.\\
\indent{\par}
{\sc Definition 5.4.}
Let $u=(u_{1}, \ldots , u_{p})\in {\cal AP}(s)$, $s\in S^{+}$.
By an {\it admissible shuffle} of $u$ we understand a sequence
of admissible partitions of $s$ of the form
$$
u=u^{0}\rightarrow u^{1} \rightarrow \ldots 
\rightarrow u^{k} \rightarrow \ldots \rightarrow
u^{p-1}=s
$$
where
$$
u^{k}=(u_{1}^{k}, \ldots , u_{n-k}^{k})=j_{l}(u^{k-1}), 
\;\;\; 1\leq l \leq n-k
$$
and
$$
j_{l}(v_{1}, \ldots , v_{m})=
(v_{1}, \ldots , v_{l}\cup v_{m}, v_{l+1}, \ldots ,v_{m-1}) , \;\;\;
1\leq l \leq m-1
$$
i.e. each transition of the shuffle amounts to moving the last subword $v_{m}$
to one of the previous subwords $v_{l}$ and then forming $v_{l}\cup v_{m}$.\\
\indent{\par}
{\it Example 3.}
The shuffle 
$$
(1,2,3,4)\rightarrow (1,2\cup 4,3) \rightarrow (1\cup 3,2\cup 4) 
\rightarrow (1\cup 2\cup 3\cup 4)
$$
is admissible if and only if $i(2,4)=i(1,3)=i(1\cup 3,2\cup 4)=1$. If
$u=(1,2,3,4)$ is admissible with
$u_{k}=s_{k}$ for all $k$, then clearly $i(2,4)=i(1,3)=1$, so we are left with
the condition
$i(1\cup 3,2\cup 4)=1$, which implies that we must have
$\delta_{2}=\delta_{3}=1$, i.e $s_{2}=s_{3}=z$. In turn, the shuffle
$$
(1,2,3,4)\rightarrow (1,2,3\cup 4) \rightarrow (1\cup 3\cup 4,2) 
\rightarrow (1\cup 2\cup 3\cup 4)
$$
is admissible if and only if $i(1\cup 3 \cup 4,2)=1$ 
which implies that we must have $\delta_{2}=1$, i.e. $s_{2}=z$.\\
\indent{\par}
{\sc Definition 5.5.}
Let $a(u)$ denote the number of admissible shuffles of u,
where $u\in {\cal AP}(s)$ (in that case $ 1\leq a(u)\leq (n-1)!$). 
If $u\notin {\cal AP}(s)$, we set $a(u)=0$.\\
\indent{\par}
{\sc Proposition 5.6.}
{\it Let $ u=(u_{1}, \ldots , u_{p})\in {\cal AP}(s)$, $p\geq 1$. Then}
\begin{eqnarray*}
a(u)&=&i(u) \;\; {\rm if}\;\; p\leq 2\\
a(u)&=&
\sum_{k=1}^{p-1}i(u_{k},u_{p}) 
a(u_{1}, \ldots , u_{k}\cup u_{p}, \ldots , u_{p-1})\;\; {\rm if}\;\;
p>2
\end{eqnarray*}
{\it Proof.}
This recurrence formula is an easy consequence of Definitions 5.4-5.5.
\hfill $\Box$.\\
\indent{\par}
{\it Example 4.}
For simplicity, assume that $u$ consists of one-letter words.
We get
$$
a(1,2,3)=i(1,3)i(1\cup 3,2)+i(2,3)i(1,2\cup 3)=
\delta_{2}+1
$$
which can be seen to agree with $m(1,2,3)$ (cf. Example 1). In a similar manner,
we get
\begin{eqnarray*}
a(1,2,3,4)&=&i(1,4)i(1\cup 4,3)i(1\cup 3\cup 4,2) + i(1,4)i(2,3)i(1\cup 4,2\cup 3)\\
&+& i(2,4)i(1,3)i(1\cup 4,2\cup 3) + i(2,4)i(2\cup 4,3)i(1,2\cup 3\cup 4)\\
&+& i(3,4)i(1,3\cup 4)i(1\cup 3\cup 4,2) + i(3,4)i(2,3\cup 4)i(1,2\cup 3\cup 4)\\
&=& \delta_{2}\delta_{3}+\delta_{2}\delta_{3}+\delta_{2}\delta_{3}+\delta_{3}
+\delta_{2}+1\\
&=&
\left\{
\begin{array}{lll}
1 & {\rm if} & s=s_{1}w^{2}s_{4}\\
2 & {\rm if} & s\in \{s_{1}wzs_{4}, s_{1}zws_{4}\} \\
6 & {\rm if} & s=s_{1}z^{2}s_{4}
\end{array}
\right.
\end{eqnarray*}
where $s_{1},s_{4}\in \{z,w\}$. \\
\indent{\par}
In order to find a connection between the 
M\"{o}bius function $m(v)$ in terms of the number of admissible
shuffles of $v$, $v\in {\cal AP}(s)$, we first need to express
$m(v)$ in terms of the M\"{o}bius functions of words $s'$, where
$l(s')<l(s)$.\\
\indent{\par}
{\sc Proposition 5.7.}
{\it If $v\in {\cal AP}(s)$ and $v<s$, where $s\in S^{+}$, then}
$$
m(v)=-\sum_{v\leq u <s}m(v|u)
$$
{\it where}
$$
m(v|u)=m(v^{1}|u_{1})\ldots m(v^{p}|u_{p})
$$
{\it for $u=(u_{1}, \ldots , u_{p})$ and $v^{k}=v\cap u_{k}$ is the partition
of $u_{k}$ consisting of those words of $v$ whose union gives $u_{k}$.}\\[5pt]
{\it Proof.}
Applying formula (5.8) to every $L(u_{k})$ 
on the RHS of (5.4), we get
$$
L(u)=\sum_{v^{1}\leq u_{1}}\ldots \sum_{v^{p}\leq u_{p}}
m(v^{1}|u_{1})\ldots m(v^{p}|u_{p})M(v^{1})\ldots M(v^{p})
$$
which, in view of (5.2), gives
$$
m(v|u)=m(v^{1}|u_{1})\ldots m(v^{p}|u_{p})
$$
using arguments similar to those in the proof of Proposition 5.2.
Therefore, we get
\begin{eqnarray*}
L(s) 
&=& 
M(s)-\sum_{0\leq u <s}L(u)\\
&=&
M(s) - 
\sum_{0\leq u <s}\sum_{v\leq u}m(v|u)M(v)\\
&=& 
M(s) - 
\sum_{0\leq v <s}\sum_{v\leq u <s}m(v|u)M(v)
\end{eqnarray*}
which gives the desired formula for $m(v)$.
\hfill
$\Box$\\
\indent{\par}
{\sc Definition 5.8.}
We say that $u$ {\it covers} $v$, where $u,v\in P$, 
$u\leq v$, if the segment $[v,u]$ contains two elements.
An {\it atom} in $P$ is an element that covers 
$0_{s}$ (a minimal element of $P$) for some $s$, 
and a {\it dual atom} is an element that is covered by 
$1_{s}$ (a maximal element of $P$) for some $s\in S^{+}$ (see [R1]). 
Denote by $D(s)$ the set of dual atoms covered by $s$.\\
\indent{\par}
{\sc Proposition 5.9.}
{\it If $v\in {\cal AP}(s)$ and $v<s$, where $s\in S^{+}$, then} 
\begin{equation}
\label{5.10}
m(v)=-\sum_{\stackrel{u\in D(s)}
{\scriptscriptstyle v_{1}, v_{2}\;{\rm separated}}}
m(v|u)
\end{equation}
{\it where the summation runs over all dual atoms $u=(u_{1},u_{2})$
of $s$ in which $v_{1}\subset u_{1}$, $v_{2}\subset u_{2}$, i.e.
the first two subwords of $v$ are separated in $u$.}\\[5pt]
{\it Proof.}
Clearly, (5.10) holds for $b(v)=2$ since in that case $LHS=-i(v)=-1$
and $RHS=m(v_{1}|v_{1})m(v_{2}|v_{2})=-1\cdot 1=-1$. Suppose (5.10)
holds for every $v\in {\cal AP}(s)$ with $2\leq b(v)\leq p-1$. We will
show that it then holds for $b(v)=p$. Let 
$v=(v_{1}, \ldots , v_{p})\in {\cal AP}(s)$. We need to show that
$$
E:=m(v)+\sum_{\stackrel{u\in D(s)}
{\scriptscriptstyle v_{1}, v_{2}\;{\rm separated}}}
m(v|u)=0.
$$
We claim that 
\begin{equation}
\label{5.11}
E=-\sum_{\stackrel{u\geq v, b(u)=r}
{\scriptscriptstyle v_{1},v_{2}\subset u_{1}}}
m(v|u)-
\sum_{\stackrel{u\geq v}
{\scriptscriptstyle b(u)>r}}
m(v|u)
\end{equation}
for $r=2, \ldots , p-1$. In view of Proposition 5.7, we have
$$
E=-\sum_{\stackrel{u\geq v, b(u)=2}
{\scriptscriptstyle v_{1},v_{2}\subset u_{1}}}
m(v|u)-
\sum_{\stackrel{u\geq v}
{\scriptscriptstyle b(u)>2}}
m(v|u),
$$
thus (5.11) holds for $r=2$.
Suppose now that (5.11) holds for $r=k$, use multiplicativity 
of $m(v|u)$ in the first sum, namely
$$
m(v|u)=m(v^{1}|u_{1})\ldots m(v^{k}|u_{k})
$$
and apply the hypothesis that (5.10) holds
for $b(v)\leq p-1$ to $m(v^{1}|u_{1})$ (note that $b(v^{1})\leq p-1$)
to get
\begin{eqnarray*}
E&=&-\sum_{\stackrel{u\geq v, b(u)=k}
{\scriptscriptstyle v_{1},v_{2}\subset u_{1}}}
\sum_{\stackrel{t\in D(u_{1})}
{\scriptscriptstyle v_{1},v_{2}\;{\rm separated}}}
m(v^{1}|t)m(v^{2}|u_{2})\ldots m(v^{k}|u_{k})
-\sum_{\stackrel{u\geq v}
{\scriptscriptstyle b(u)>k}}
m(v|u)\\
&=&
-\sum_{\stackrel{u\geq v, b(u)=k+1}
{\scriptscriptstyle v_{1},v_{2}\subset u_{1}}}
m(v|u)
-
\sum_{\stackrel{u\geq v}
{\scriptscriptstyle b(u)>k+1}}
m(v|u).
\end{eqnarray*}
This equality is justified as follows. If $u\geq v$, $b(u)=k$ and 
$v_{1},v_{2}\subset u_{1}$ and $t=(u_{1}',u_{1}'')$ 
is a dual atom of $u_{1}$ which separates $v_{1} $ and $v_{2}$, 
then the partition $(u_{1}',u_{1}'',u_{2}, \ldots , u_{k})$ is admissible.
In turn, if $u \geq v$, $b(u)=k+1$, then either $v_{1},v_{2}\subset u_{1}$,
or $v_{1}\subset u_{1}$ and $v_{2}\subset u_{2}$.
In the second case, there exists exactly one pair $(u',t)$, where
$u'=(u_{1}',\ldots , u_{k}')\in {\cal AP}(s)$ and a dual atom 
$t$ of the first subword $u_{1}'$ such that the resulting partition
is $u$ (take $u_{1}$ and $u_{2}$, form $u_{1}\cup u_{2}$, 
the remaining $u_{l}$'s keep the same and let the dual atom 
of $u_{1}\cup u_{2}$ be $(u_{1},u_{2})$).
Therefore, (5.11) holds for all $2\leq r \leq p-1$. But, if $r=p-1$, then
(5.11) takes the form
$$
E=-m(v_{1},v_{2}|v_{1}\cup v_{2})m(v_{3}|v_{3})\ldots m(v_{p}|v_{p})-
m(v_{1}|v_{1})\ldots m(v_{p}|v_{p})=0
$$
which finishes the proof.\hfill $\Box$\\
\indent{\par}
{\sc Corollary 5.10.}
{\it For all $v\in {\cal AP}(s)$, where $s\in S^{+}$,
we have}
\begin{equation}
\label{5.12}
m(v)=(-1)^{b(v)-1}a(v)
\end{equation}
{\it Proof.}
Clearly, $m(s)=1=a(s)$ and if $v=(v_{1},v_{2})$, then 
$m(v)=-i(v)=-a(v)$.
We will show that if (5.12) holds for $b(v) \leq p-1$, then 
it holds for $b(v)=p$. Let $b(v)=p$ and use Proposition 5.9 and 
the inductive assumption to get
\begin{eqnarray*}
m(v)&=&
-
\sum_{\stackrel{v=v^{1}\cup v^{2}}
{\scriptscriptstyle v_{1}\in v^{1},v_{2}\in v^{2}}}
(-1)^{b(v^{1})+b(v^{2})-2}i(u_{1},u_{2})a(v^{1})a(v^{2})\\
&=&
(-1)^{b(v)-1}
\sum_{\stackrel{v=v^{1}\cup v^{2}}
{\scriptscriptstyle v_{1}\in v^{1},v_{2}\in v^{2}}}
i(u_{1},u_{2})a(v^{1})a(v^{2})
\end{eqnarray*}
where $(u_{1},u_{2})$'s are the dual atoms which appear on the RHS 
of (5.10).
Now, from Definition 5.4 we can see that in order to count
all admissible shuffles of  
$v=(v_{1}, \ldots , v_{p}) \in {\cal AP}(s)$
it is enough to count all admissible shuffles 
of $v$ which lead to the dual atom $(u_{1},u_{2})$  
after $p-2$ transitions such that $v_{1}\subset u_{1}$
and $v_{2}\subset u_{2}$ (here $u_{1},u_{2}$ correspond to
$v_{1}^{p-2}, v_{2}^{p-2}$ of Definition 5.4) --
this computation gives the product $a(v^{1})a(v^{2})$
-- and then take into account only those such pairs which give
an admissible partition of $s$ (this gives $i(u_{1},u_{2})$).
\hfill $\Box$\\
\indent{\par}
{\it Example 5.}
We apply Corollary 5.10 to compute $m(u)$ needed for
the cumulants of lowest order.
For simplicity, we write 
$L(i_{1}\ldots i_{n})=L(s_{i_{1}}\ldots s_{i_{n}})$ 
and 
$M(i_{1}\ldots i_{n})=M(s_{i_{1}}\ldots s_{i_{n}})$. 
We obtain
\begin{eqnarray*}
L(1)
&=&
M(1)\\
L(12)
&=&
M(12)-M(1)M(2)\\
L(123)
&=&
M(123)-M(12)M(3)-M(1)M(23)-\delta_{2}M(13)M(2)\\
& +& (1+\delta_{2})M(1)M(2)M(3)\\
L(1234)
&=&
M(1234)-M(123)M(4)-\delta_{2}M(134)M(2)-\delta_{3}M(124)M(3)\\
&-&
M(1)M(234)- M(12)M(34) -\delta_{2}\delta_{3}M(14)M(23)\\
&-&
\delta_{2}\delta_{3}M(13)M(24) +(1+\delta_{3})M(12)M(3)M(4)\\
&+&
\delta_{2}(1+\delta_{3})M(13)M(2)M(4)
+2\delta_{2}\delta_{3}M(14)M(2)M(3)\\
&+&(1+\delta_{2}\delta_{3})M(1)M(23)M(4)
+\delta_{3}(1+\delta_{2})M(1)M(24)M(3)\\
&+&(1+\delta_{2})M(1)M(2)M(34)
-(1+\delta_{2}+\delta_{3}+3\delta_{2}\delta_{3})M(1)M(2)M(3)M(4).
\end{eqnarray*}
One can recognize some coefficients computed in Examples 1,2,4.
Also note that if all $\delta$'s are equal to $1$
(i.e. $s=s_{1}z^{2}s_{4}$), we get classical cumulants, whereas
if all $\delta$'s are equal to $0$ (i.e. $s=s_{1}w^{2}s_{4}$), we get
boolean cumulants.\\
\myownsection
\begin{center}
{\sc 6. Semigroup algebras}
\end{center}
Let us now briefly recall
some basic facts on the free semigroup algebra 
${\cal A}(S)$ and the Banach algebra $l^{1}(S,W)$, where 
$W$ is a weight function. For more on this subject
see [P].

The {\it free semigroup algebra} of $S$, denoted ${\cal A}(S)$,
is the set of functions
$$
f:\; S\rightarrow {\bf C}
$$
equipped with the usual addition and convolution multiplication given by
\begin{equation}
\label{6.1}
f\star g (s)=\sum_{uv=s}f(u)g(v).
\end{equation}
Note that, in general, $f\star g \neq g \star f$, so the algebra 
${\cal A}(S)$ is noncommutative.

It is convenient to identify elements $f$ of both algebras with 
formal sums
$$
\sum_{s\in S}f(s)s
$$
with multiplication
\begin{equation}
\label{6.2}
\sum_{t\in S}f(t)t \; \sum_{r\in S}g(r)r=\sum_{s\in S}\sum_{tr=s}f(t)g(r)s
\end{equation}
The algebra ${\cal A}(S)$ is our noncommutative analog of the 
algebra of formal power series ${\bf C}[[z]]$ in the variable $z$.
The unit of the algebra is denoted ${\bf 1}$, where
${\bf 1}(s)=1$ if $s=1$ and otherwise is equal to zero.\\
\indent{\par}
{\sc Proposition 6.1.}
{\it Let $f\in {\cal A}(S)$ and assume that $f(1)=1$. Then 
$f$ is invertible in ${\cal A}(S)$, i.e. there exists
$f^{-1}\in {\cal A}(S)$, such that $f^{-1}\star f =f\star f^{-1}={\bf 1}$,
$f^{-1}(1)=1$ and}
$$
f^{-1}(s)=\sum_{p=1}^{l(s)}(-1)^{p}
\sum_{s=u_{1}\ldots u_{p}}f(u_{1})\ldots f(u_{p})
$$
{\it for any $s\in S^{+}$, where the second sum is taken over all 
factorizations of $s$.}\\
\indent{\par}
We omit the proof since it is a straightforward computation.\\
\indent{\par}
{\sc Definition 6.2.}
For given $f\in {\cal A}(S)$ and $g\in {\cal A}(S(z))$, we define
$f_{\star g}$ to be the function from 
${\cal A}(S)$ which agrees with $f$ on $S(z)$ and, on $S\setminus S(z)$,
is given by
$$
f_{\star g} (z_{1}w_{1}z_{2}w_{2}\ldots z_{p-1}w_{p-1}z_{p})
$$
$$
= \sum_{u_{1}v_{1}=z_{1}}\ldots \sum_{u_{p-1}v_{p-1}=z_{p-1}}
f(u_{1}w_{1}u_{2}w_{2}\ldots u_{p-1}w_{p-1}z_{p})
g(v_{1})\ldots g(v_{p-1})
$$
for $w_{1},\ldots , w_{p-1}\in S(w)\setminus \{1\}$,
$z_{1}, z_{p}\in S(z)$, $z_{2}, \ldots , z_{p-1}\in S(z)\setminus \{1\}$,
$p\geq 2$.\\
\indent{\par}
One can think of $f_{\star g}$ as a ``composition'' of $g$ and $f$ in which
$w_{1}, \ldots , w_{p-1}$ are replaced by $g(z)w_{1}, \ldots , g(z)w_{p-1}$. 
Thus, we can informally write
$$
f_{\star g}=\sum_{s\in S}f(s)s_{\star g}
$$
where $s_{\star g}$ agrees with $s$ for $s\in S(z)$ and
$$
s_{\star g}=
z_{1}g(z)w_{1}z_{2}g(z)w_{2}\ldots z_{p-1}g(z)
w_{p-1}z_{p}
$$
for 
$$
s=z_{1}w_{1}z_{2}w_{2}\ldots z_{p-1}w_{p-1}z_{p}\in S\setminus S(z),
$$
with the assumptions as in Definition 6.2.

Note also that if $f\in {\cal A}(S), g\in {\cal A}(S(z))$, with $g(1)\neq 0$, 
the following implication holds:
\begin{equation}
\label{6.3}
f_{\star g}=h \;\Rightarrow \; f=h_{\star g^{-1}},
\end{equation}
which is a straightforward consequence of Definition 6.2.

By a {\it weight function} on $S$ we understand a (real-valued) positive 
function on $S$ which is submultiplicative, i.e.
$$
W(st)\leq W(s)W(t) \;\; \forall \;s,t\in S
$$
and by $l^{1}(S,W)$ we denote the Banach space of all functions 
$f: S\rightarrow {\bf C}$ which are finite with respect to the norm
$$
\parallel f \parallel_{W} =\sum_{s\in S}W(s)|f(s)|,
$$
and which becomes a Banach algebra under the convolution multiplication (6.1),
see [P]. The $l^{1}$-{\it semigroup algebra} of $S$ [B], denoted 
$l^{1}(S)$, is obtained if $W(s)=1$ for all $s\in S$.\\
\indent{\par}
{\sc Proposition 6.3.}
{\it Let $f\in l^{1}(S,W)$ with $f(1)=1$ and $g\in l^{1}(S(z),W)$, where
$W(1)=1$ and $W(st)=W(ts)$ for all $t,s\in S$, and let $Q>1/2$.
Then the following implications hold:\\
(i) if $\parallel f \parallel _{W}\leq 2-1/Q$, 
then $\parallel f^{-1}\parallel _{W} \leq Q$,\\
(ii) if $\parallel g \parallel _{W} <Q$, then 
$\parallel f_{\star g} \parallel _{\widetilde{W}} 
\leq \parallel f \parallel _{W}$,\\
where $\widetilde{W}(s)=W(s)Q^{-m(s)}$ and $m(s)$ is the number of 
$w$'s in the word $s$.}\\[5pt]
{\it Proof.}
Using Proposition 6.1, triangle inequality and submultiplicativity of $W$,
we arrive at
\begin{eqnarray*}
\parallel f^{-1}\parallel _{W} &=&1+\sum_{s\neq 1}W(s)
|\sum_{p=1}^{l(s)}(-1)^{p}\sum_{s=u_{1}\ldots u_{p}}f(u_{1})\ldots u_{p})|\\
&\leq& 1+ \sum_{s\neq 1}W(s)|f(s)| +
\sum_{s\neq 1}W(s)|\sum_{s=u_{1}u_{2}}f(u_{1})
f(u_{2})| + \ldots\\
&=& \parallel f \parallel _{W} +(\parallel f \parallel _{W}-1)^{2} + \ldots\\
&=& \frac{1}{2-\parallel f \parallel _{W}}
\end{eqnarray*}
from which (i) follows. Now,
$$
\parallel f_{\star g}\parallel _{W}
=
\sum_{s\in S(z)}W(s)|f(s)|+\sum_{s\in S\setminus S(z)}W(s)|f_{\star g}(s)|,
$$
and an estimate of the second sum is needed.
In the sums below we will always assume that all $w_{k}$'s
belong to $S(w)\setminus \{1\}$ and that all $z_{k}$'s belong
to $S(z)$ (additional restrictions on $z_{k}$'s will be given
explicitly) without further mention. Therefore
\begin{eqnarray*}
\sum_{s\in S\setminus S(z)}\widetilde{W}(s)|f_{\star g}(s)|& = &
\sum_{p=2}^{\infty}
\sum_{w_{1}, \ldots , w_{p-1}}
\sum_{\stackrel{z_{1}, \ldots , z_{p}}
{\scriptscriptstyle {z_{2}, \ldots , z_{p-1} \neq 1}}}
\widetilde{W}(z_{1}w_{1}z_{2}\ldots z_{p-1}w_{p-1}z_{p})\\
&&\times 
|f_{\star g}(z_{1}w_{1}z_{2}\ldots z_{p-1}w_{p-1}z_{p}|\\
&\leq &
\sum_{p=2}^{\infty}
\sum_{w_{1}, \ldots , w_{p-1}}
\sum_{\stackrel{z_{1}, \ldots , z_{p}}
{\scriptscriptstyle {z_{2}, \ldots , z_{p-1} \neq 1}}}
\sum_{u_{1}v_{1}=z_{1}}\ldots \sum_{u_{p-1}v_{p-1}=z_{p-1}}\\
&& 
\widetilde{W}(u_{1}w_{1}u_{2}\ldots u_{p-1}w_{p-1}z_{p})
|f(u_{1}w_{1}u_{2}\ldots u_{p-1}w_{p-1}z_{p}|\\
&& \times 
\prod_{l=1}^{p-1}W(v_{l})|g(v_{l})|\\
&\leq &
\sum_{r\in S\setminus S(z)}\widetilde{W}(r)
|f(r)|\parallel g \parallel _{W}^{m(r)}
\end{eqnarray*}
where we used submultiplicativity of $W$ and $W(st)=W(ts)$. When 
we use the definition of $\widetilde{W}$ and combine the above estimate
with the sum over $S(z)$, we get (ii).
\hfill $\Box$\\
\myownsection
\begin{center}
{\sc 7. Moment and cumulant generating functions}
\end{center}
In this section we introduce moment and cumulant generating functions
associated with filtered convolution and derive a connection between them.

Let us first establish a moment-cumulant formula,
which expresses the moments in terms of the cumulants {\it and}
the moments (in contrast to the inversion formula which expresses
the cumulants in terms of the moments only). 
This formula turns out very useful in deriving
an explicit form of the cumulant generating function.

For given $s\in S\setminus S(z)$ we denote by
$C_{0}(s)$ the subset of cumulant subwords of $s$ 
which contain the first letter $w$ in $s$. In turn, if
$r\in C_{0}(s)$, by $W_{r}(s)$ we will denote the set of subwords 
of $s\setminus r$ of maximal lenght which lie between the $w$-legs of $r$
or before the first $w$ of $r$. Note that they have to be 
$z$-words since otherwise $r$ would have an inner $w$, which is not possible
since $r$ is a cumulant word.\\
\indent{\par}
{\sc Lemma 7.1.} ({\sc Moment-cumulant formula})
{\it For each $s\in S\setminus S(z)$ we have the following
moment-cumulant formula:} 
\begin{equation}
\label{7.1.}
M(s)=\sum_{r\in C_{0}(s)}L(r)
\prod_{v\in W_{r}(s)}
M(v) 
\times
M(s\setminus (r\cup \bigcup_{v\in W_{r}(s)}v)),
\end{equation}
{\it Proof.}
From Definition 4.1 we have
$$
M(s)=\sum_{p=1}^{l(s)}\sum_{u=(u_{1}, \ldots , u_{p})\in {\cal AP}(s)}
L(u_{1})\ldots L(u_{p})
$$
for $s\in S\setminus S(z)$. Since $s$ must contain at least one $w$, 
in each of the summands on the RHS of the above formula we must have
one cumulant, say $L(u_{j})$, such that $u_{j}$ contains the first
$w$ of the word $s$. Denote, for fixed $u$, this $u_{j}$ by $r=r(u)$. 
From the definition of admissible partitions of $s$,
letters from $s\setminus r$ lying between 
two $w$-legs of $r$ cannot be connected 
with letters lying
between another pair of $w$-legs of $r$ since otherwise they would be separated
by a $w$. 
The same is true for the letters lying to the left of the first $w$.
Moreover, these letters have to be $z$'s.
Therefore, the product of cumulants corresponding to such  
$z$-words $v\in W_{r}(s)$ 
have to be taken over all subpartitions of $v$, giving
$M(v)$. The same applies to the word $v\in W_{r}(s)$ formed from all
$z$'s which are to the left of the first $w$ and are not in $r$.
Altogether, these products of cumulants give
$$
\prod_{v\in W_{r}(s)}M(v).
$$
The remaining cumulants from the product $L(u_{1})\ldots L(u_{p})$
involve only letters lying to the right of the last $w$ of the word $r$
which are not in $r$. The product of them gives
$$
M(s\setminus (r\cup \bigcup_{u\in W_{r}(s)}v))
$$
which completes the proof.\hfill $\Box$\\
\indent{\par}
{\it Example 1.}
To illustrate the moment-cumulant formula, let us give a diagram
corresponding to one of the summands on the RHS of (7.1).
We choose the word $s=z^{3}wz^{3}wz^{2}wz$ (long enough to 
show some general features of the combinatorics involved).
The diagram\\ 
\unitlength=1mm
\special{em:linewidth 0.4pt}
\linethickness{0.4pt}
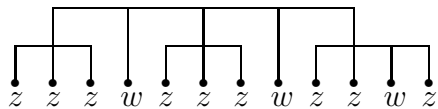
\begin{picture}(120.00,25.00)(-37.00,5.00)
\put(10.00,15.00){\line(1,0){10.00}}
\put(30.00,15.00){\line(1,0){10.00}}
\put(50.00,15.00){\line(1,0){15.00}}
\put(15.00,20.00){\line(1,0){40.00}}
\put(10.00,15.00){\line(0,-1){5.00}}
\put(20.00,15.00){\line(0,-1){5.00}}
\put(30.00,15.00){\line(0,-1){5.00}}
\put(40.00,15.00){\line(0,-1){5.00}}
\put(50.00,15.00){\line(0,-1){5.00}}
\put(60.00,15.00){\line(0,-1){5.00}}
\put(65.00,15.00){\line(0,-1){5.00}}
\put(15.00,20.00){\line(0,-1){10.00}}
\put(25.00,20.00){\line(0,-1){10.00}}
\put(35.00,20.00){\line(0,-1){10.00}}
\put(45.00,20.00){\line(0,-1){10.00}}
\put(55.00,20.00){\line(0,-1){10.00}}
\put(10.00,10.00){\circle*{1.00}}
\put(15.00,10.00){\circle*{1.00}}
\put(20.00,10.00){\circle*{1.00}}
\put(25.00,10.00){\circle*{1.00}}
\put(30.00,10.00){\circle*{1.00}}
\put(35.00,10.00){\circle*{1.00}}
\put(40.00,10.00){\circle*{1.00}}
\put(45.00,10.00){\circle*{1.00}}
\put(50.00,10.00){\circle*{1.00}}
\put(55.00,10.00){\circle*{1.00}}
\put(60.00,10.00){\circle*{1.00}}
\put(65.00,10.00){\circle*{1.00}}
\put(9.00,7.00){$z$}
\put(14.00,7.00){$z$}
\put(19.00,7.00){$z$}
\put(24.00,7.00){$w$}
\put(29.00,7.00){$z$}
\put(34.00,7.00){$z$}
\put(39.00,7.00){$z$}
\put(44.00,7.00){$w$}
\put(49.00,7.00){$z$}
\put(54.00,7.00){$z$}
\put(59.00,7.00){$w$}
\put(64.00,7.00){$z$}
\end{picture}
$\;$\\
corresponds to
$$
L(zwzwz)\times M(z^{2})M(z^{2}) \times M(zwz)
$$
and the upper line connects all letters associated with the cumulant,
whereas the lower line connects all letters associated with the 
moments. \\
\indent{\par}
{\it Remark 1.}
If $s\in S(z)$ (the case not treated in Proposition 7.1), we get the classical 
moment-cumulant formula
\begin{equation}
\label{7.2}
M(z^{n})=\sum_{k=1}^{n}{n-1 \choose k-1}L(z^{k})M(z^{n-k}),
\end{equation}
where $n\geq 1$, i.e. $M(z^{n})$, $L(z^{n})$ 
are classical moments and cumulants of order $n$, 
respectively. Now, formulas (7.1) and (7.2)
cover the cases of all $s\in S$. In (7.2) we already counted 
the number of ways, namely ${n-1 \choose k-1}$, 
of choosing all subwords $r$ of $s=z^{n}$ which are equal to $z^{k}$ as words
and contain the first $z$ in $s$  --
these subwords correspond to the cumulant $L(z^{k})$ -- 
and the product of the remaining cumulants on the RHS of (4.1)
gives $M(z^{n-k})$ since all partitions of $s\setminus r$
are allowed.\\
\indent{\par}
{\it Remark 2.}
In turn, the boolean moment-cumulant formula [Sp-W] is a special case of 
(7.1) and
\begin{equation}
\label{7.3}
M(w^{n})=\sum_{k=1}^{n}L(w^{k})M(w^{n-k}),
\end{equation}
where $M(w^{n})$ and $L(w^{n})$ are boolean moments and cumulants
of order $n$, respectively. In this case, all cumulant subwords 
of $s=w^{n}$ containing the first $w$ are of the form 
$r=s_{1}\ldots s_{k}=w^{k}$ (for each $k$ there is only one such subword),
the product over $W_{r}(s)$ in (7.1) disappears and the remaining
moment is therefore equal to $M(s\setminus r)=M(w^{n-k})$.\\
\indent{\par}
Note that there is a formal similarity between our formula and 
the moment-cumulant formula in the conditionally-free case
[Bo-Le-Sp]. Nevertheless, there is a substantial difference
between the two cases -- our formula involves not only 
non-crossing partitions, which later gives rise to some classical
features in the generating functions, namely they are  
analogs of {\it exponential} generating functions.\\
\indent{\par}
{\sc Definition 7.2.}
Let $(M(s))_{s\in S}$ and $(L(s))_{s\in S}$
be the moments and cumulants associated with the state  $\widehat{\phi}$.
The corresponding moment and cumulant
generating functions are defined to be the elements
of the algebra ${\cal A}(S)$ given by the formal sums
\begin{eqnarray*}
M\{z,w\}&=&\sum_{s\in S}\frac{M(s)}{n(s)!}s,\\
L\{z,w\}&=&\sum_{1\neq s\in S}\frac{L(s)}{n(s)!}s,
\end{eqnarray*}
respectively, where 
\begin{equation}
\label{7.4}
n(s)!=n_{1}!n_{2}!\ldots n_{p}!\;\; {\rm for}\;\;
s=z^{n_{1}}w^{k_{1}}z^{n_{2}}w^{k_{2}}\ldots w^{k_{p-1}}z^{n_{p}}
\end{equation}
with $n_{1},n_{p}\in {\bf N}_{0}$ and 
$k_{1},n_{2}, k_{2}, \ldots , n_{p-1},k_{p-1}\in {\bf N}$.\\
\indent{\par}
Now, let us use Definition 6.2 to introduce new notations
\begin{eqnarray}
\label{7.5}
L^{\star} \{z,w\}&=&L_{\star M}\{z,w\}\\
\label{7.6}
M_{\star} \{z,w\}&=&M_{\star M^{-1}}\{z,w\}
\end{eqnarray}
where we take $f_{\star g}$ with $f=L\{z,w\}$
and $g=M\{z,0\}$ in (7.5) and 
with $f=M\{z,w\}$ and $g= M^{-1}\{z,0\}$
in (7.6). Here,
$$
M\{z,0\}=\sum_{s\in S(z)}\frac{M(s)}{l(s)!}s
$$
i.e. $M\{0,z\}$ is the restriction of $M\{z,w\}$ to the 
support $S(z)$ (then $n(s)=l(s)$) and can be treated as a formal
power series in $z$ representing the classcial moment generating
function. In a similar way we define $L\{z,0\}$, the classical
cumulant generating function,
as well as $M\{0,w\}$ and $L\{0,w\}$ (in these two cases,
by restricting the support to $S(w)$).
Note that $M^{-1}\{z,0\}$, the inverse of $M\{z,0\}$, 
exists since $M(1)=1$. Also note that $M_{\star}\{z, 0\}=M\{z, 0\}$
and $L^{\star}\{z,0\}=L\{z, 0\}$.

Moreover, for $f=f\{z,w\}\in {\cal A}(S)$,
we will use a special notation for the difference
\begin{equation}
\label{7.7}
\delta f\{z,w\}=f\{z,w\}- f\{z, 0\},
\end{equation}
representing the ``deviation from the classical case''
and apply this notation to 
$\delta M\{z,w\}$, $\delta M_{\star} \{z,w\}$,
$\delta L\{z,w\}$ and $\delta L^{\star} \{z,w\}$.

Using these notations we can write down formulas which connect
the cumulant generating function with the moment generating function.\\
\indent{\par}
{\sc Theorem 7.3.}
{\it The moment and cumulant generating functions of Definition 7.2
satisfy the relation}
\begin{equation}
\label{7.8}
\delta M\{z,w\}=
\delta L^{\star}\{z,w \}
M\{ z, w\}
\end{equation}
{\it with the multiplication of formal sums given by (6.2).}\\[5pt]
{\it Proof.}
Note that
$$
\delta M\{z,w\}=
\sum_{s\in S\setminus S(z)}\frac{M(s)}{n(s)!}s
$$
and thus let us consider the RHS of the above equation.
Use the moment-cumulant formula of Lemma 7.1
for each $M(s)$ with $s\in S\setminus S(z)$ (there is at least
one $w$ in each $s$).
On the RHS of (7.1) we have to compute the number of ways 
in which the same {\it word} is obtained
by taking different {\it subwords} of the word $s$ 
(i.e. different subsequences). 
Suppose $s$ is of the form
given by (7.4). If $r\in C_{0}(s)$, then $r$ either ends with a $z$, 
i.e. is of the form
\begin{equation}
\label{7.9}
r=z^{i_1}w^{k_1}z^{i_2}w^{k_2}\ldots z^{i_{m-1}}w^{k_{m-1}}z^{i_{m}}, \;\;
{\rm with}\;\; m\leq p ,\; 1\leq i_{m}\leq n_{m}
\end{equation}
or ends with a $w$, i.e. is of the form
\begin{equation}
\label{7.10}
r=z^{i_1}w^{k_1}z^{i_2}w^{k_2}\ldots z^{i_{m}}w^{q_{m}}, \;\;
{\rm with}\;\; m\leq p-1 ,\; 1\leq q_{m}\leq k_{m},
\end{equation}
where all powers are assumed to be positive except perhaps $i_{1}$, which may
be equal to zero. Note that $r$ has to assume one of these forms 
(i.e. the powers of $w$ have to coincide with those in $s$ up to some place)
since otherwise the partition $u=(u_{1}, \ldots , u_{p})$ corresponding to the
product of cumulants $L(u_{1})\ldots L(u_{p})$ would not be 
admissible.

We will use the multiindex notation
$$
{n(s)\choose n(r)}=
{n_{1} \choose i_{1}}{n_{2}\choose i_{2}}\ldots {n_{p}\choose i_{p}}
$$
where 
$$
n(s)=(n_{1},n_{2}, \ldots , n_{p}),\;\;\;
n(r)=(i_{1}, i_{2}, \ldots , i_{p})
$$
for $s$ of the form (7.4) and $r$ given by (7.9) or (7.10),
where we set $i_{m+1}=\ldots = i_{p}=0$.

Thus, if we want to include in the summation only those 
words corresponding to $r\in C_{0}(s)$ which are distinct, we get
$$
M(s)=\sum_
{\stackrel{r\in C_{0}(s)}
{\scriptscriptstyle {\rm distinct\; words}}}
L(r){n(s)\choose n(r)}
\prod_{v\in W_{r}(s)}M(v)\; M(s\setminus (r\cup \bigcup_{v\in W_{r}(s)}v))
$$
which leads to the equation
\begin{eqnarray*}
\delta M\{z,w\}
&=&
\sum_{s\in S\setminus S(z)}
\sum_{\stackrel{r\in C_{0}(s)}
{\scriptscriptstyle {\rm distinct \;words}}}
\frac{L(r)}{n(r)!}
\frac{1}{(n(s)-n(r))!}
\prod_{v\in W_{r}(s)}M(v)
\\
&\times &
M(s\setminus (r\cup \bigcup_{v\in W_{r}(s)}v))s
\end{eqnarray*}
where
$$
(n(s)-n(r))!=(n_{1}-i_{1})!(n_{2}-i_{2})!\ldots (n_{p}-i_{p})!.
$$
In order to demonstrate (7.8), we need to show that
$$
\delta M\{z,w\} (s)=\delta L^{\star}\{z,w\} M\{z,w\} (s)
$$
(treated as elements of ${\cal A}(S)$) for every $s\in S\setminus S(z)$. 

Let us write $\delta L^{\star}\{z,w\}$ informally as
$$
\delta L^{\star}=
\sum_{r\in S\setminus S(z)}\frac{L(r)}{n(r)!}r^{\star}
$$
where
$$
r^{\star}=
z^{i_{1}}M\{z,0\}w^{k_{1}}z^{i_{2}}M\{z,0\}w^{k_{2}}\ldots
z^{i_{m-1}}M\{z,0\}w^{k_{m-1}}z^{i_{m}}
$$
for $r$ of the form (7.9) (we allow $i_{m}=0$ which means that
the case (7.10) is also covered in this notation). Using this notation, we can 
write
$$
\delta L^{\star}\{z,w\} M\{z,w\}=
\sum_{r\in S\setminus S(z)}\sum_{t\in S}
\frac{L(r)}{n(r)!}\frac{M(t)}{n(t)!}r^{\star}t= f_{1}\{z,w\}+f_{2}\{z,w\}
$$
where, in the last equality, we have 
split the sum on the RHS of the above formula into two sums:
the first one in which between the last $w$ from $r^{\star}$ and 
the first $w$ from $t$ there is a $z$ and the second one,
in which between the last $w$ from $r^{\star}$ and the first $w$
from $t$ there are no $z$'s.

We have
$$
f_{1}\{z,w\}(z^{n_{1}}w^{k_{1}}z^{n_{2}}w^{k_{2}}\ldots
z^{n_{p-1}}w^{k_{p-1}}z^{n_{p}})=
$$
\begin{eqnarray*}
&=&\sum_{m=1}^{p}\sum_{i_{1}=0}^{n_{1}}\sum_{i_{2}=1}^{n_{2}} \ldots 
\sum_{i_{m-1}=1}^{n_{m-1}} \sum_{i_{m}=1}^{n_{m}}
\frac{1}{i_{1}!\ldots i_{m}!}\frac{1}{(n_{1}-i_{1})!\ldots (n_{m-1}-i_{m-1})!}\\
&\times &
L(z^{i_{1}}w^{k_{1}}z^{i_{2}}w^{k_{2}}\ldots z^{i_{m-1}}w^{k_{m-1}}z^{i_{m}})
M(z^{n_{1}-i_{1}})\ldots M(z^{n_{m-1}-i_{m-1}})\\
&\times&
\frac{1}{(n_{m}-i_{m})!n_{m+1}!\ldots n_{p}!}
M(z^{n_{m}-i_{m}}w^{k_{m}}z^{n_{m+1}}\ldots z^{n_{p}})
\end{eqnarray*}
and
$$
f_{2}\{z,w\}(z^{n_{1}}w^{k_{1}}z^{n_{2}}w^{k_{2}}\ldots
z^{n_{p-1}}w^{k_{p-1}}z^{n_{p}})=
$$
\begin{eqnarray*}
&=&\sum_{m=1}^{p-1}\sum_{i_{1}=0}^{n_{1}}\sum_{i_{2}=1}^{n_{2}} \ldots 
\sum_{i_{m-1}=1}^{n_{m-1}} \sum_{q_{m-1}=1}^{k_{m-1}}
\frac{1}{i_{1}!\ldots i_{m-1}!}\frac{1}{(n_{1}-i_{1})!
\ldots (n_{m-1}-i_{m-1})!}\\
&\times& 
L(z^{i_{1}}w^{k_{1}}z^{i_{2}}w^{k_{2}}\ldots z^{i_{m-1}}w^{q_{m-1}})
M(z^{n_{1}-i_{1}})\ldots M(z^{n_{m-1}-i_{m-1}})\\
&\times& \frac{1}{n_{m}!n_{m+1}!\ldots n_{p}!}
M(w^{k_{m-1}-q_{m-1}}z^{n_{m}}w^{k_{m}}\ldots z^{n_{p}}).
\end{eqnarray*}
It is not hard to see that 
\begin{eqnarray*}
f_{1}\{z,w\}(s)&=&
\sum_{\stackrel{r\in C_{0}'(s)}
{\scriptscriptstyle {\rm distinct\;words}}}
\frac{L(r)}{n(r)!}\frac{1}
{(n(s)-n(r))!}
\prod_{v\in W_{r}(s)}M(v)\\
&\times& M(s\setminus (r\cup \bigcup_{v\in W_{r}(s)}v))
\end{eqnarray*}
and
\begin{eqnarray*}
f_{2}\{z,w\}(s)&=&
\sum_{\stackrel{r\in C_{0}''(s)}
{\scriptscriptstyle {\rm distinct \; words}}}
\frac{L(r)}{n(r)!}\frac{1}
{(n(s)-n(r))!}\prod_{v\in W_{r}(s)}M(v)\\
&\times& M(s\setminus (r\cup \bigcup_{v\in W_{r}(s)}v))
\end{eqnarray*}
where $C_{0}'(s)$ and $C_{0}''(s)$ are 
the subsets of $C_{0}(s)$ which consist
of those subwords which end with a $z$ and $w$, respectively
(suumations are taken only over subwords which give distinct words).

Thus
$$
f_{1}\{z,w\}(s)+f_{2}\{z,w\}(s)=\delta M (s)
$$
which finishes the proof. 
\hfill $\Box$\\
\indent{\par}
{\sc Corollary 7.4.}
{\it The cumulant generating function takes the form}
\begin{equation}
\label{7.11}
L\{z,w\}=L\{z,0\}+
(M_{\star}\left\{z,w \right\}- M\left\{z,0\right\})
M_{\star}^{-1}\left\{ z, w\right\}
\end{equation}
{\it where the notation (7.5)-(7.7) is used.}\\[5pt]
{\it Proof.}
By multiplying (7.8) from the right by the inverse of $M\{z,w\}$,
which exists since $M(1)=1$, we get
$$
\delta M \{z,w\} M^{-1}\{z,w\}= \delta L^{\star} \{z,w\}
$$
but now, using (6.3) we can get rid of $\star$ on the RHS of this equation.
This leads to
$$
\delta M_{\star} \{z,w\} M_{\star}^{-1}\{z,w\}= \delta L \{z,w\}
$$
which is equivalent to (7.11).
\hfill $\Box$\\
\indent{\par}
{\it Remark.}
The classical moment and cumulant generating functions 
can be identified with $M\{z, 0\}$ (the Fourier transform) 
and $L\{z,0\}$ (its logarithm), respectively,
and are related through the classical equation
$$
L\{z, 0\}={\rm log}M\{z,0\}
$$
which can be derived from (7.2) in the usual manner. In turn, the boolean
moment and cumulant generating functions can be identified with $M\{0,w\}$
and $L\{0,w\}$, respectively, and are related through the equation
$$
L\{0,w\}=(M\{0,w\}-1)M^{-1}\{0,w\}
$$
which is a special case of (7.11) since $M_{\star}\{0,w\}=
M\{0,w\}$, $L\{0,0\}=0$ and $M_{\star}^{-1}\{0,w\}=M\{0,w\}^{-1}$.\\
\indent{\par}
{\sc Corollary 7.5.}
{\it If $\parallel M\{z,w\}\parallel _{W}<2-1/Q$ and
$L\{z,0\}\in l^{1}(S(z),W)$, where
$W$ and $Q$ satisfy the assumptions of Proposition 6.3., then
equation (7.11) holds in $l^{1}(S,\widetilde{W})$.}\\[5pt]
{\it Proof.}
Let us first note that for any state $\widehat{\phi}$,
there exists $W$ such that the assumptions of this Corollary
are satisfied for $M_{\widehat{\phi}}\{z,w\}$ and 
$L_{\widehat{\phi}}\{z,0\}$. Now, in view of Proposition 6.3(i),
these assumptions imply that
$$
\parallel M^{-1}\{z,w\}\parallel _{W}\leq Q
$$
which, by Proposition 6.3(ii), gives
$$
\parallel M_{\star}^{-1}\{z,w\}\parallel_{\widetilde{W}}\leq 
\parallel M^{-1}\{z,w\}\parallel \leq Q.
$$
Besides, $\parallel M\{z,0\}\parallel _{W}< 2-1/Q$ implies
that
$$
\parallel M_{\star}\{z,w\}\parallel _{\widetilde{W}}
\leq \parallel M\{z,w\}\parallel _{W}<2-1/Q
$$
and therefore, the RHS of (7.11) is an element of $l^{1}(S,\widetilde{W})$
and thus (7.11) holds in $l^{1}(S,\widetilde{W})$.
\hfill $\Box$\\
\indent{\par}
Note that the weight function $W$ plays a role similar to 
the radius of convergence of a power series. 
Thus, if we say that the formula for the cumulant generating 
function, derived on the level of semigroup algebra ${\cal A}(S)$,
namely (7.11), holds on the ``analytic'' level, i.e. there exists
$W$ such that (7.11) holds in $l^{1}(S,W)$, it is analogous to saying
that for a generating function in the form
of a formal power series from ${\bf C}[[z]]$, there exists $R>0$
such that this power series becomes convergent in the circle of radius 
$R$. \\
\begin{center}
{\sc REFERENCES}
\end{center}
[A] {\sc M.~Anshelevich}, ``Partition-dependent stochastic
measures and $q$-deformed cumulants'', MSRI Preprint No.
2001-021, Berkeley, 2001.\\[3pt]
[B] {\sc B.~A.~Barnes, J.~Duncan}, ``The Banach algebra $l^{1}(S)$'',
{\it J.~Funct.~Anal.} {\bf 18} (1975), 96-113.\\[3pt]
[Bo-Le-Sp] {\sc M.~Bozejko, M.~Leinert, R.~Speicher},
``Convolution and limit theorems for con\-di\-tio\-na\-lly 
free random variables'',
{\it Pacific J. Math.} {\bf 175} (1966), 357-388.\\[3pt]
[G-S] {\sc A.~B.~Ghorbal, M.Sch\"{u}rmann}, ``On the algebraic formulation 
of non-com\-mu\-ta\-tive pro\-ba\-bi\-li\-ty theory'', 
preprint, Universite Henri Poincare, 1999.\\[3pt]
[F] {\sc U.~Franz}, ``Unification of boolean, monotone, anti-monotone
and tensor independence and Levy processes'', preprint 4/2001, 
Ernst-Moritz-Arndt-Universitat Greifswald, 2001.\\[3pt]
[F-L] {\sc U.~Franz, R.~Lenczewski},
``Limit theorems for the hierarchy of freeness'', 
{\it Prob. Math. Stat.} {\bf 19} (1999), 23-41.\\[3pt]
[L1] {\sc R.~Lenczewski}, ``Unification of independence in
quantum probability'', {\it Inf.~Dim. Anal.~Quant.~Probab. \& Rel.~Top.}
{\bf 1} (1998), 383-405.\\[3pt]
[L2] {\sc R.~Lenczewski}, ``Filtered stochastic calculus'',
{\it Inf. Dim.~Anal. Quant.~Probab. \& Rel.~Top.}, to appear.\\[3pt]
[L3] {\sc R.~Lenczewski}, ``Filtered random variables, bialgebras and 
convolutions'', {\it J. Math. Phys.}, to appear.\\[3pt]
[Lo] {\sc M.~Lothaire}, {\it Combinatorics on Words}, Encyclopedia
of Mathematics and its Applications, Addison-Wesley, London, 1983.\\[3pt]
[Mu1] {\sc N.~Muraki}, ``Monotonic independence, monotonic central limit theorem
and monotonic law of large numbers'',{\it Inf.~Dim.~Anal.~
Quant.~Probab.~Rel.~Top.},to appear.\\[3pt]
[Mu2] {\sc N.~Muraki}, ``Monotonic convolution and monotonic Levy-Hincin
formula'', preprint, 2000.\\[3pt]
[N] {\sc A.~Nica}, ``A one-parameter family of transforms, linearlizing
convolution laws for probability distributions'', {\it Comm.~Math.~Phys}
{\bf 168} (1995), 187-207.\\[3pt]
[P] {\sc T.~Palmer}, {\it Banach Algebras and the General Theory of 
*-Algebras, Vol.I}, Cambridge University Press, 1994.\\[3pt]
[R1] {\sc G-C.~Rota}, ``On the foundations of combinatorial theory I.
Theory of M\"{o}bius functions'', {\it Z.~Wahrscheinlichkeitstheorie
Verw. Gebiete} {\bf 2} (1964), 340-368.\\[3pt]
[R2] {\sc G-C.~Rota}, {\it Finite Operator Calculus},
Academic Press, New York, 1975.\\[3pt]
[S] {\sc M.~Sch\"{u}rmann}, ``Direct sums of tensor products and
non-commutative independence'', {\it J. Funct. Anal.} {\bf 133}
(1995), 1-9.\\[3pt]
[Sp] {\sc R.~Speicher}, ``Multiplicative functions on the lattice
of non-crossing partitions and free convolution'',
{\it Math.~Ann.} {\bf 298} (1994), 611-628.\\[3pt]
[Sp-W] {\sc R.~Speicher, R.~Woroudi}, ``Boolean convolution'', {\it
Fields Institute Commun.} {\bf 12} (1997), 267-279.\\[3pt]
[V1] {\sc D.~Voiculescu}, ``Symmetries of some reduced free product 
${\cal C}^{*}$-algebras'', in {\it Operator Algebras and their 
Connections with Topology and Ergodic Theory}, Lecture
Notes in Math. 1132, Springer, Berlin, 1985, 556-588.\\[3pt]
[V2] {\sc D.~Voiculescu}, ``Addition of certain non-commuting random
variables'', {\it J.~Funct. Anal.} {\bf 66} (1986), 323-346.\\[3pt]
\end{document}